\documentclass[onecolumn]{autart}    
\usepackage{graphicx}          
\usepackage{subfigure}                            
\usepackage{mathtools}                               
\usepackage{amsmath}
\usepackage{mathrsfs}
\usepackage{amssymb}
\usepackage{bm}
\usepackage{float}
\usepackage{hyperref}
\usepackage{color}
\newcommand{\mcl}[1]{\mathcal{#1}}

\newcommand{\pr}{\hbox{\sf P}}

\newcommand{\ep}{\hbox{\sf E}}
\newcommand{\var}{\hbox{\sf Var}\,}
\newcommand{\cov}{\hbox{\sf Cov}\,}

\newcommand{\id}{\hbox{\bf 1}}

\newtheorem{problem}{Problem}
\newtheorem{remark}{Remark}
\allowdisplaybreaks[4]
\begin{document}

\begin{frontmatter}

\title{Mean-variance portfolio selection with dynamic attention behavior
in a hidden Markov model} \vspace{-3em}


\thanks[footnoteinfo]{Corresponding author Jiaqin~Wei.}
\author[*]{Yu Zhang}\ead{2997179078@qq.com},
\author[Z. Jin]{Zhuo Jin}\ead{zhuo.jin@mq.edu.au},
\author[*]{Jiaqin Wei}\ead{jqwei@stat.ecnu.edu.cn},
\author[G. Yin]{George Yin}\ead{gyin@uconn.edu}\vspace{-0.5em}

\address[*]{Key Laboratory of Advanced Theory and Application in Statistics and
Data Science-MOE, School of Statistics, East China Normal University, Shanghai 200062, China}\vspace{-0.5em}
\address[Z. Jin]{Department of Actuarial Studies and Business Analytics, Macquarie University, 2109, NSW, Australia}\vspace{-0.5em}
\address[G. Yin]{Department of Mathematics, University of Connecticut, Storrs, CT 06269-1009, USA}\vspace{-1.5em}

\begin{keyword}                           
Mean-variance; Dynamic attention behavior; Extended HJB equation; Markov chain approximation; Hidden Markov model.              
\end{keyword}                             

\begin{abstract}                          
In this paper, we study closed-loop equilibrium strategies for mean-variance
portfolio selection problem in a hidden Markov model with dynamic
attention behavior. In addition to the investment strategy, the investor's
attention to news is introduced as a control of the accuracy of the
news signal process. The objective is to find equilibrium strategies by numerically
solving an extended HJB equation by using Markov chain approximation
method. An iterative algorithm is constructed and its convergence is established.
Numerical examples are also provided to illustrate the results.\vspace{-1.5em}
\end{abstract}

\end{frontmatter}

\section{Introduction}
Since the introduction of mean-variance (MV, for short) portfolio selection problem
by \cite{Mh:52}, much progress has been made. Various extensions have been investigated including the dynamic
MV problem in multi-period model and continuous-time models;
see, e.g., \cite{LdNw:00} and \cite{ZxLd:00}. Certain
limitations of a single diffusion and the existence of `regimes' that
switch among themselves in the underlying market are realized, a system commonly
referred to as the regime-switching (RS, for short) model was used to discuss MV problems (see, e.g., \cite{ZxYg:03}). Usually,
some states of the market cannot be observed by an investor. In the financial market with partially observable states, hidden Markov chain is often
used to describe the evolution of the unobservable market states. The portfolio
optimization problem in the hidden Markov RS financial
market has been investigated by \cite{ErStBa:10} and \cite{YzYgZq:15}.

It is well-known that the dynamic MV optimization problem
is time-inconsistent in the sense that Bellman's optimality principle
does not hold. In \cite{Sr:55}, where a time-inconsistent problem
within a game theoretic framework was studied, the author viewed the
time-inconsistent optimization problem as a multi-person game and
looked for a subgame perfect Nash equilibrium point. The precise definition
of the equilibrium concept in continuous time within the class of
closed-loop strategies was provided for the first time in \cite{EiLa:06}
and \cite{EiPt:08}. They investigated the optimal investment and
consumption problem with hyperbolic discounting in continuous-time
deterministic and stochastic models, respectively. Recently, there
has been increasing interests in finding the so-called equilibrium
strategies for the time-inconsistent optimization problems; see \cite{BtMa:10}, \cite{HyJhZx:12}, \cite{WtJzWj:19}, \cite{WjWt:17}, among others.

In most of the existing literature on hidden Markov model, it is assumed
that the accuracy of signal collection is fixed in the financial market,
and investors acquire information passively in the sense that they
do not control the quality of the information they collect. However,
it is reasonable to introduce costly information-acquisition problems,
which can help investors obtain accurate information, but at the expense
of decreasing her current wealth. In \cite{DjKr:87},
the author examined the demand for information and derived the equilibrium
price of information. An investor who faces uncertainty about stock-return predictability was considered in \cite{Xy:01}. The authors solved the portfolio-choice and information-acquisition problem. The relation between costly information acquisition and the excess covariance of asset prices was studied in \cite{Vl:06}. In \cite{HlLh:07} where a dynamic portfolio-choice problem with
static costly information choice was studied, the author assumed that
the accuracy and frequency of information was chosen at one time only.
However, in \cite{GaRa:18}, an investor's attention was assumed time varying. Recently, a dynamic portfolio-choice
problem with dynamic information acquisition was considered in \cite{AdHm:20}. At each point in time, the investor optimally chose the quantity of information that she needed, which results in a dynamic trade-off problem between asset and attention allocation.

In this paper, we consider the time-inconsistent MV portfolio
selection in a hidden Markov model with dynamic attention behavior. Following
the concept of investor's attention to news in \cite{AdHm:20}, we
introduce the dynamic attention behavior viewed as an endogenous
control into the hidden Markov model. We assume that the investor
knows more accurate information about the current state of the market
as her attention to news increases. Under this framework, the
extended Hamilton-Jacobi-Bellman (HJB, for short) equation
cannot be solved in closed-form.
So our main objective is to find numerical approximation of the resulting
control problem. To obtain the numerical solutions, we adopt the Markov chain
approximation method developed in \cite{KhDp:13},
which has been used to find optimal controls
in various complex stochastic systems in finance, insurance, and other fields, see \cite{JzYhYg:13}, \cite{JzYgWf:13}, \cite{JzYgZc:12}, \cite{Ss:08}, \cite{SqYgZz:06}, \cite{TkYg:16}, among others.
To the best of our knowledge,  this paper is the first work to study extended HJB equation numerically
by using the Markov chain approximation method, while the aforementioned references are all about solving the classical HJB equation. We build an iterative algorithm and obtain the convergence of the algorithm. Besides, numerical examples are provided to illustrate the results.
The numerical results show that an investor with more wealth tends to acquire more information, which is consistent with those in \cite{AdHm:20}.

The remainder of this paper is organized as follows. Section 2 introduces
the model with the investor's attention to news. The Markov chain
approximation method is considered in Section 3. In section 4, we
consider the approximation of optimal controls and establish the convergence
of the algorithm. Section 5 gives one numerical example for illustration.
Finally, Section 6 concludes the paper with further remarks.

\section{Model}
Let $T>0$ be a fixed finite time horizon and $\boldsymbol{W}_{1}(\cdot)=\left(W_{1}^{1}(\cdot),\dots,W_{1}^{d}(\cdot)\right)^{\top}$
be a $d$-dimensional standard Brownian motion, where the symbol
$\top$ indicates the transpose.\footnote{In the following, we shall suppress the suffix $(\cdot)$ and the time variable for processes and functions whenever there is no 
 confusion.}
Let $\alpha$ be a continuous-time Markov chain independent
of $\boldsymbol{W}_{1}$ valued in the finite set $\mathscr{M}=\{1,2,\dots,m\}$
and denote by $\boldsymbol{Q}=(q^{ij})_{m\times m}$ the generator
of $\alpha$. Let $(\Omega,\mcl F,\pr)$
be a complete probability space on which $\boldsymbol{W}_{1}$,
$\alpha$, and all the random variables and processes in the
rest of the paper are defined.

We consider a financial market consisting of one bond and $d$ stocks
within the time horizon $[0,T]$. The price of the risk-free bond
$S_{0}$ satisfies
\begin{equation*}
\setlength{\abovedisplayskip}{3pt}
\setlength{\belowdisplayskip}{3pt}
dS_{0}=r\left(s,\alpha\right)S_{0}ds,\quad s\in(0,T],
\end{equation*}
where $S_{0}(0)=s_{0}>0$ and $r(\cdot,i)\geq0$ is the risk-free return rate in state $i$ for each
$i=1,2,\dots,m$.

For $l=1,2,\ldots,d$, the price of the $l$-th stock $S_{l}$
is given by
\begin{equation*}
\setlength{\abovedisplayskip}{3pt}
\setlength{\belowdisplayskip}{3pt}
dS_{l} =\mu_{l}(s,\alpha)S_{l}ds+\sum_{j=1}^{d}\sigma_{lj}(s,\alpha)S_{l}dW_{1}^{l},\quad s\in(0,T],
\end{equation*}
where $S_{l}(0)=s_{l}>0$, $\mu_{l}(\cdot,i)$ and $\sigma_{lj}(\cdot,i),i=1,2,\dots,m$
are the expected return 
 and volatility rate of the $l$th risky
asset corresponding to state $i$, respectively.

In our framework, instead of having full information of the Markov
chain $\alpha$, the investor can only observe it in white
noise.
However, we allow the investor to have opportunity to actively
learn about regime predictability.

Similar to \cite{AdHm:20}, we assume that the investor may acquire
$n(s)$ signals of equal precision $Z^{j}(s),j=1,\ldots,n(s)$
at time $s$, and each signal represents information from market, which is given by
\begin{equation*}
\setlength{\abovedisplayskip}{3pt}
\setlength{\belowdisplayskip}{3pt}
dZ^{j}=\zeta(\alpha)ds+\sigma dB^{j},\quad j=1,\ldots,n(s),
\end{equation*}
where $B^{j}$ for $j=1,\dots,n(s)$ are independent real-valued standard Brownian motions that are independent of $\boldsymbol{W}_{1}$. Here, the function $\zeta$ is given by a vector $\boldsymbol{\zeta}=(\zeta_{1}, \zeta_{2}, ..., \zeta_{m})'$, so that $\zeta(\alpha)=\left\langle \alpha, \boldsymbol{\zeta} \right\rangle$ where $\left\langle \cdot,\cdot \right\rangle$ denotes the scalar product in $\mathbb{R}^{m}$.
Then there is a standard real-valued Brownian motion $B$ so that the aggregated signal $Z({s})\coloneqq\sum_{j=1}^{n(s)}dZ^{j}(s)/n(s)$ acquired by the investor satisfies
\begin{equation}
\setlength{\abovedisplayskip}{3pt}
\setlength{\belowdisplayskip}{3pt}
dZ=\zeta(\alpha)ds+\sigma/\sqrt{n(s)}dB.\label{eq:2-4}
\end{equation}
Therefore the investor can choose the number of signals
$n(s)$ to control the accuracy of the aggregated
signal. If we set $\pi (s)= n(s)/ \sigma^2$. Then clearly, the diffusion coefficient in (\ref{eq:2-4}) is $1/\sqrt {\pi(s)}$.
To take this idea further, we shall not restrict $\pi(s)$ to take only
discrete values when $n(s)$ is integer-valued, but consider it to take all real values in an arbitrary interval.
This can be considered as a ``fluid approximation.'' It is more convenient for us to work with a real-valued control function than an integer-valued function.

Thus, we choose diffusion term as a control and consider the following signal process
\begin{equation*}
\setlength{\abovedisplayskip}{3pt}
\setlength{\belowdisplayskip}{3pt}
dZ=\zeta\left(\alpha \right)ds+1/\sqrt{\pi}dB,\quad s\in(0,T],
\end{equation*}
with $Z(0)=0$. Here, we assume that $\pi\in [\epsilon,M]$, where $\epsilon,M>0$ are two arbitrary constants. The boundedness of $\pi$ is consistent with \cite{Sc:03}, which argued that investors have limited information-processing capacity.

We further assume that the investor can get accurate signals
at the cost of decreasing her current wealth, i.e., the investor
faces a dynamic trade-off problem of asset and attention allocation. Denote by $\bar{K}$ the total information cost and assume that $\bar{K}=K(\pi)X$, where the per-unit-of-wealth cost function $K(\pi)$ is increasing and convex in attention 
\cite{AdHm:20}.
Let $u_l$ be the dollar amount invested in the $l$th stock. Then the wealth of the investor, denoted by $X$, evolves as
\begin{equation}
\setlength{\abovedisplayskip}{3pt}
\setlength{\belowdisplayskip}{3pt}
dX = b(s,X,\alpha,\boldsymbol{u},\pi)ds+\boldsymbol{\sigma}(s,X,\alpha,\boldsymbol{u})d\boldsymbol{W}_{1},
\label{eq:2-6}
\end{equation}
where $X(0)=x_{0}$, $\boldsymbol{u}\coloneqq\left(u_{1},\ldots,u_{d}\right)^{\top}$,
\[
\setlength{\abovedisplayskip}{3pt}
\setlength{\belowdisplayskip}{3pt}
\begin{aligned}
b\left(s,X,\alpha,\boldsymbol{u},\pi\right) & \coloneqq r\left(s,\alpha\right)X+\boldsymbol{\theta}^{\top}\left(s,\alpha\right)\boldsymbol{u}-\bar{K},\\
\boldsymbol{\sigma}\left(s,X,\alpha,\boldsymbol{u}\right) & \coloneqq\boldsymbol{u}^{\top}\left(\sigma_{lj}\left(s,\alpha\right)\right)_{1\leq l,j\leq d}.
\end{aligned}
\]
Here, $\boldsymbol{\theta}(\cdot,i)\coloneqq\left(\mu_{1}(\cdot,i)-r(\cdot,i),\ldots,\mu_{d}(\cdot,i)-r(\cdot,i)\right)^{\top}$, for $i=1,2,\dots,m$.

Next, let us recall some results about Wonham's filter. To proceed,
we denote by $I_{E}$ the indicator function of the event $E$ and define $p^{i}(t)\coloneqq I_{\{\alpha(t)=i\}}$, and
\begin{equation*}
\setlength{\abovedisplayskip}{3pt}
\setlength{\belowdisplayskip}{3pt}
\varphi^{i}(t)\coloneqq\pr\left(\alpha(t)=i\mid Z(s),0\leq s\leq t\right).
\end{equation*}
Following \cite{TkYg:16}, we know that $\varphi^{i}(t),i=1,2,\dots,m$ is the probability conditioned
on the observation $\sigma\left(Z(s),0\leq s\leq t\right)$. Since $\varphi^{i}(t)\geq0$
and $\sum_{i=1}^{m}\varphi^{i}(t)=1$, it
is sufficient to work with $\boldsymbol{\varphi}\coloneqq\left(\varphi^{1},\dots,\varphi^{m-1}\right)^{\top}$. Denote by
\begin{equation*}
\setlength{\abovedisplayskip}{3pt}
\setlength{\belowdisplayskip}{3pt}
S_{m-1} \coloneqq\left\{ \boldsymbol{\varphi}=\left(\varphi^{1},\dots,\varphi^{m-1}\right)^{\top}:\varphi^{i}\geq0,\right.
 \left.\sum_{i=1}^{m-1}\varphi^{i}\leq1\right\},
\end{equation*}
and write $\varphi^{m}=1-\sum_{i=1}^{m-1}\varphi^{i}$ for $\boldsymbol{\varphi}\in S_{m-1}$.

It was shown in \cite{Ww:64} that the posterior probability
$\boldsymbol{\varphi}$
satisfies the following system of stochastic differential equation
\begin{equation*}
\setlength{\abovedisplayskip}{3pt}
\setlength{\belowdisplayskip}{3pt}
d\boldsymbol{\varphi}=\tilde{\boldsymbol{b}}\left(\boldsymbol{\varphi}\right)ds+\tilde{\boldsymbol{\sigma}}\left(\boldsymbol{\varphi},\pi\right)d\hat{W}_{2},\quad s \in (0,T],
\end{equation*}
where $\boldsymbol{\varphi}(0)=\boldsymbol{\varphi}_{0}\coloneqq\left(\varphi_{0}^{1},\dots,\varphi_{0}^{m-1}\right)^{\top}$ is the initial distribution of $\alpha$,
\[
\setlength{\abovedisplayskip}{3pt}
\setlength{\belowdisplayskip}{3pt}
\begin{aligned}\tilde{\boldsymbol{b}}\left(\boldsymbol{\varphi}\right) & \coloneqq\left(\sum_{j=1}^{m}q^{j1}\varphi^{j},\dots,\sum_{j=1}^{m}q^{j(m-1)}\varphi^{j}\right)^{\top},\\
\tilde{\boldsymbol{\sigma}}\left(\boldsymbol{\varphi},\pi\right) & \coloneqq\left(\sqrt{\pi}\varphi^{1}\left(\zeta(1)-\bar{\zeta}\left(\boldsymbol{\varphi}\right)\right),\dots,\sqrt{\pi}\varphi^{m-1}\left(\zeta(m-1)-\bar{\zeta}(\boldsymbol{\varphi})\right)\right),
\end{aligned}
\]
$\bar{\zeta}\left(\boldsymbol{\varphi}\right)\coloneqq\sum_{i=1}^{m}\zeta(i)\varphi^{i}$, and the innovation process $\hat{W}_{2}$, independent of $\boldsymbol{W}_{1}$, is defined as
\begin{equation*}
\setlength{\abovedisplayskip}{3pt}
\setlength{\belowdisplayskip}{3pt}
\hat{W}_{2}(t)\coloneqq\int_{0}^{t}\sqrt{\pi(s)}dZ(s)-\int_{0}^{t}\sqrt{\pi(s)}\bar{\zeta}\left(\boldsymbol{\varphi}(s)\right)ds.	
\end{equation*}
\qquad With $p^{i}(t)$, (\ref{eq:2-6}) can be written as
\begin{equation*}
\setlength{\abovedisplayskip}{3pt}
\setlength{\belowdisplayskip}{3pt}
dX =\sum_{i=1}^{m}p^{i}\left[b(s,X,i,\boldsymbol{u},\pi)ds+\boldsymbol{\sigma}(s,X,i,\boldsymbol{u})d\boldsymbol{W}_{1}\right].
\end{equation*}
By replacing the hidden state $p^{i}$ by its estimate $\varphi^{i}$, we obtain an estimated process of the original
process. For notational simplicity, we still denote by $X$
the estimated process. Letting $\boldsymbol{Y}\coloneqq\left(X,\boldsymbol{\varphi}\right)^{\top}$, then we can get a completely observed system with the following compacted form
\begin{equation}
\setlength{\abovedisplayskip}{3pt}
\setlength{\belowdisplayskip}{3pt}
d\boldsymbol{Y}=\boldsymbol{f}(s,\boldsymbol{Y},\boldsymbol{u},\pi)ds+\boldsymbol{\varSigma}(s,\boldsymbol{Y},\boldsymbol{u},\pi)d\hat{\boldsymbol{v}},\label{eq:2-10}
\end{equation}
where $\boldsymbol{Y}(0)=\boldsymbol{y}_{0}\coloneqq(x_{0},\boldsymbol{\varphi}_{0})^{\top}$, $\hat{\boldsymbol{v}}\coloneqq\left(\boldsymbol{W}_{1},\hat{W}_{2}\right)^{\top}$,
\[
\setlength{\abovedisplayskip}{3pt}
\setlength{\belowdisplayskip}{3pt}
\boldsymbol{f}(s,\boldsymbol{Y},\boldsymbol{u},\pi)\coloneqq\left(\bar{b}(s,\boldsymbol{Y},\boldsymbol{u},\pi),\tilde{\boldsymbol{b}}(\boldsymbol{\varphi})\right)^{\top},
\]
\[
\begin{aligned}\boldsymbol{\varSigma}(s,\boldsymbol{Y},\boldsymbol{u},\pi) & \coloneqq\begin{pmatrix}\bar{\boldsymbol{\sigma}}(s,\boldsymbol{Y},\boldsymbol{u}) & 0\\
\mathbf{0} & \tilde{\boldsymbol{\sigma}}\left(\boldsymbol{\varphi},\pi\right)
\end{pmatrix}.\end{aligned}
\]
Here, we define $\bar{b}(s,\boldsymbol{Y},\boldsymbol{u},\pi)\coloneqq\sum_{i=1}^{m}b(s,x,i,\boldsymbol{u},\pi)\varphi^{i}$
and $\bar{\boldsymbol{\sigma}}(s,\boldsymbol{Y},\boldsymbol{u})\coloneqq\sum_{i=1}^{m}\boldsymbol{\sigma}(s,x,i,\boldsymbol{u})\varphi^{i}$.

Let $\mathbb{F}\coloneqq\{\mathcal{F}_{s}: s\in[0,T]\}$, where
${\mathcal F}_s=\{\boldsymbol{W}_{1}(\tilde s), Z(\tilde s), X(t): t\le \tilde s\le s\}$.
For $\boldsymbol{H}\coloneqq\mathbb{R},\mathbb{R}^{d}$, etc. and $0\leq t\leq s\leq T$, define
\[
\setlength{\abovedisplayskip}{3pt}
\setlength{\belowdisplayskip}{3pt}
\begin{gathered}
L_{\mathbb{F}}^{2}\left(t,s;\boldsymbol{H}\right)\coloneqq\left\{ \boldsymbol{X}:[t,s]\times\Omega\rightarrow\boldsymbol{H}\mid\boldsymbol{X}\text{ is }\mathbb{F}\text{-adapted,\;}\ep\left[\int_{t}^{s}|\boldsymbol{X}(v)|^{2}dv\right]<\infty\right\} ,\\
L_{\mathbb{\mathbb{F}}}^{2}\left(\Omega;C\left([t,s];\boldsymbol{H}\right)\right)\coloneqq\left\{ \boldsymbol{X}:[t,s]\times\Omega\rightarrow\boldsymbol{H}\mid\boldsymbol{X}\text{ is }\mathbb{F}\text{-adapted, }\text{has continuous paths, and\;}\ep\left[\sup_{v\in[t,s]}|\boldsymbol{X}(v)|^{2}\right]<\infty\right\} .
\end{gathered}
\]
\qquad In this paper, we restrict ourselves to closed-loop controls, that is at each $t$, $\left(\boldsymbol{u}(t),\pi(t)\right)$
is given by $\boldsymbol{u}(t) =\tilde{\boldsymbol{u}}\left(t,\boldsymbol{Y}(t)\right), \pi(t)=\tilde{\pi}\left(t,\boldsymbol{Y}(t)\right),$
where the maps $\tilde{\boldsymbol{u}}:[0,T]\times\mathbb{R}^{m}\rightarrow\mathbb{R}_{+}^{d}$
and $\tilde{\pi}:[0,T]\times\mathbb{R}^{m}\rightarrow[\epsilon,M]$ are
two Borel measurable functions. Here, we call the functions $\tilde{\boldsymbol{u}}$
and $\tilde{\pi}$ control laws.
In the following, for notational simplicity, we still write $\boldsymbol{u}$
and $\pi$ for $\tilde{\boldsymbol{u}}$ and $\tilde{\pi}$, respectively.

\begin{defn}
Admissible control laws are maps $\boldsymbol{u}$ and $\pi$ satisfying
the following conditions:
\begin{enumerate}
\item For each initial state $\left(t,\boldsymbol{y}\right)\in[0,T]\times\mathbb{R}^{m}$, (\ref{eq:2-10})
has a unique solution $\boldsymbol{Y}\in L_{\mathbb{\mathbb{F}}}^{2}\left(\Omega;C\left([t,T];\mathbb{R}^{m}\right)\right)$.
\item The processes $\boldsymbol{u}(\cdot,\boldsymbol{Y}(\cdot))\in L_{\mathbb{F}}^{2}\left(0,T;\mathbb{R}_{+}^{d}\right)$ and $\pi(\cdot,\boldsymbol{Y}(\cdot))\in L_{\mathbb{F}}^{2}\left(0,T;[\epsilon,M]\right)$.
\end{enumerate}
\end{defn}
We denote by $\boldsymbol{U}$ and $\varPi$ the class of admissible control
laws $\boldsymbol{u}$ and $\pi$, respectively.

For any initial state $\left(t,\boldsymbol{y}\right)$, the objective
of the investor is to minimize the MV cost functional
\begin{equation}
\setlength{\abovedisplayskip}{3pt}
\setlength{\belowdisplayskip}{3pt}
J\left(t,\boldsymbol{y};\boldsymbol{u},\pi\right)=\var_{t}\left[X(T)\right]-\frac{\gamma}{2}\ep_{t}\left[X(T)\right],\; t\in[0,T],\label{eq:2-11}
\end{equation}
where $\gamma>0$ is a constant and $\ep_{t}[\cdot]\coloneqq\ep[\cdot\left|\mathcal{F}_{t}\right.]$.
Letting $F(x)=x-\frac{\gamma}{2}x^{2}$ and $G(x)=\frac{\gamma}{2}x^{2}$,
we can rewrite (\ref{eq:2-11}) as
\begin{equation}
\setlength{\abovedisplayskip}{3pt}
\setlength{\belowdisplayskip}{3pt}
J\left(t,\boldsymbol{y};\boldsymbol{u},\pi\right)=\ep_{t}\left[F\left(X(T)\right)\right]-G\left[\ep_{t}\left(X(T)\right)\right].\label{eq:2-12}
\end{equation}
\qquad In summary, we need to solve the following problem

\begin{problem}
\label{prob1}At any time $t\in[0,T]$ with wealth and posterior probability
$\boldsymbol{Y}(t)=\left(X(t),\boldsymbol{\varphi}(t)\right)$, specify admissible control
laws $\boldsymbol{u}\in\boldsymbol{U}$ and $\pi\in\varPi$ that minimize
(\ref{eq:2-12}).
\end{problem}

It is well-known that Problem \ref{prob1} is time-inconsistent as
the term $G\left[\ep_{t}\left(X(T)\right)\right]$ involves a non-linear
function of the expectation. In this paper we aim to find the time-consistent
equilibrium strategies for this problem. Following \cite{BtMa:10},
we give the definition of equilibrium strategies and equilibrium value
function.

\begin{defn}
Consider admissible control laws $\bar{\boldsymbol{u}}$ and $\bar{\pi}$.
For arbitrary admissible control laws $\boldsymbol{u}\in\boldsymbol{U}$,
$\pi\in\varPi$ and a fixed real number $\varepsilon>0$, define the
control laws $\left(\boldsymbol{u}_{\varepsilon},\pi_{\varepsilon}\right)$
by
\begin{equation*}
\setlength{\abovedisplayskip}{3pt}
\setlength{\belowdisplayskip}{3pt}
\begin{pmatrix}\boldsymbol{u}_{\varepsilon}\left(s,\boldsymbol{y}\right)\\
\pi_{\varepsilon}\left(s,\boldsymbol{y}\right)
\end{pmatrix}=\begin{cases}
\begin{pmatrix}\boldsymbol{u}\left(s,\boldsymbol{y}\right)\\
\pi\left(s,\boldsymbol{y}\right)
\end{pmatrix} & \text{if}\;t\leq s<t+\varepsilon,\boldsymbol{y}\in\mathbb{R}^{m},\\
\begin{pmatrix}\bar{\boldsymbol{u}}\left(s,\boldsymbol{y}\right)\\
\bar{\pi}\left(s,\boldsymbol{y}\right)
\end{pmatrix} & \text{if}\;t+\varepsilon\leq s\leq T,\boldsymbol{y}\in\mathbb{R}^{m}.
\end{cases}
\end{equation*}
If
\begin{equation*}
\setlength{\abovedisplayskip}{3pt}
\setlength{\belowdisplayskip}{3pt}
\liminf_{\varepsilon\rightarrow0}\frac{J\left(t,\boldsymbol{y};\bar{\boldsymbol{u}},\bar{\pi}\right)-J\left(t,\boldsymbol{y};\boldsymbol{u}_{\varepsilon},\pi_{\varepsilon}\right)}{\varepsilon}\geq0
\end{equation*}
for arbitrary $\varepsilon>0$ and $\left(t,\boldsymbol{y}\right)\in[0,T]\times\mathbb{R}^{m}$,
we say that $\left(\bar{\boldsymbol{u}},\bar{\pi}\right)$ is an equilibrium
control law. The corresponding equilibrium value function $V$ is
given by $V\left(t,\boldsymbol{y}\right)\coloneqq J\left(t,\boldsymbol{y};\bar{\boldsymbol{u}},\bar{\pi}\right)$.
\end{defn}
\qquad For any twice continuously differentiable function $\phi(\cdot,\cdot,\cdot):\mathbb{R}_{+}\times\mathbb{R}\times S_{m-1}\mapsto\mathbb{R}$,
we define
\[
\begin{aligned}
\mathscr{L}^{\boldsymbol{u},\pi}\phi(t,\boldsymbol{y}) & \coloneqq\frac{\partial\phi}{\partial t}+\frac{\partial\phi}{\partial x}\bar{b}\left(t,\boldsymbol{y},\boldsymbol{u},\pi\right)+\frac{1}{2}\frac{\partial^{2}\phi}{\partial x^{2}}\bar{\boldsymbol{\sigma}}(t,\boldsymbol{y},\boldsymbol{u})\bar{\boldsymbol{\sigma}}^{\top}(t,\boldsymbol{y},\boldsymbol{u})\\
 & +\sum_{i=1}^{m-1}\frac{\partial\phi}{\partial\varphi^{i}}\sum_{j=1}^{m}q^{ji}\varphi^{j}+\frac{1}{2}\pi\sum_{i,k=1}^{m-1}\varphi^{i}\varphi^{k}\left(\zeta(i)-\bar{\zeta}(\boldsymbol{\varphi})\right)\left(\zeta(k)-\bar{\zeta}(\boldsymbol{\varphi})\right)\frac{\partial^{2}\phi}{\partial\varphi^{i}\partial\varphi^{k}}.
\end{aligned}
\]
From the general theory of Markovian time-inconsistent stochastic
control (see, e.g., \cite{BtMa:10}), we know that the equilibrium
value function $V\left(t,\boldsymbol{y}\right)$ satisfies the following
system of extended HJB equation
\begin{equation}
\setlength{\abovedisplayskip}{3pt}
\setlength{\belowdisplayskip}{3pt}
\begin{aligned}\underset{\boldsymbol{u},\pi}{\inf}\left\{ -\frac{\gamma}{2}\frac{\partial^{2}g(t,\boldsymbol{y})}{\partial x^{2}}\bar{\boldsymbol{\sigma}}(t,\boldsymbol{y},\boldsymbol{u})\bar{\boldsymbol{\sigma}}^{\top}(t,\boldsymbol{y},\boldsymbol{u})+\mathscr{L}^{\boldsymbol{u},\pi}V(t,\boldsymbol{y})-\frac{\gamma}{2}\pi\sum_{i,k=1}^{m-1}\varphi^{i}\left(\zeta(i)-\bar{\zeta}(\boldsymbol{\varphi})\right)\varphi^{k}\left(\zeta(k)-\bar{\zeta}(\boldsymbol{\varphi})\right)\frac{\partial^{2}g(t,\boldsymbol{y})}{\partial\varphi^{i}\partial\varphi^{k}}\right\}  & =0\\
V\left(T,\boldsymbol{Y}(T)\right) & =x,\label{eq:2-15}\\
\mathscr{L}^{\bar{\boldsymbol{u}},\bar{\pi}}g(t,\boldsymbol{y}) & =0,\\
g\left(T,\boldsymbol{Y}(T)\right) & =x,
\end{aligned}
\end{equation}
where $\bar{\boldsymbol{u}}$ and $\bar{\pi}$ are control laws which
realize the infimum in the first equation of (\ref{eq:2-15}).

\begin{remark}
Although we adopt the similar idea of dynamic attention behavior as in \cite{AdHm:20}, we work under a hidden Markov model and use the Wonham filter to estimate
the posterior probability. This is different from \cite{AdHm:20} where the author uses the Kalman filter.
As we know, the Wonham filter is a non-linear filter, which makes it more difficult to handle.
Besides, we study equilibrium strategies for a MV problem via extended HJB equation while \cite{AdHm:20} considers a classical utility maximization problem.
\end{remark}

\section{Discrete-time approximation scheme }

In this section, we present a numerical algorithm to solve (\ref{eq:2-15})
using Markov chain approximation method. First, we introduce a pair
of step sizes $h=(h_{1},h_{2})$. Here, $h_{1}>0$ is the discretization
parameter for state variables, $h_{2}>0$ is the step size for time
variable such that $N_{h_{2}}=T/h_{2}$ is an integer without loss of generality.
Define $\boldsymbol{S}_{h_{1}}\coloneqq\left\{ (k^{1}h_{1},\dots,k^{m}h_{1})^{\top}:k^{i}=0,\pm1,\dots,i=1,\dots,m\right\} $
and let $\{\boldsymbol{\xi}_{n}^{h},n<\infty\}$ be a discrete-time
controlled Markov chain with state space $\boldsymbol{S}_{h_{1}}$.
Let $\boldsymbol{u}^{h}=(\boldsymbol{u}_{0}^{h},
\boldsymbol{u}_{1}^{h},\dots)$
and $\boldsymbol{\pi}^{h}=
(\pi_{0}^{h},\pi_{1}^{h},\dots)$
denote the sequences of $\mathbb{R}_{+}^{d}$ and $[\epsilon,M]$ valued
random variables that are the control actions at time $0,1,\dots$,
respectively.

Denote by $\pr^{h}\left((\boldsymbol{y},\boldsymbol{z})\left|\boldsymbol{r},c\right.\right)$
the probability that $\boldsymbol{\xi}$ transits from state $\boldsymbol{y}$
at time $nh_{2}$ to state $\boldsymbol{z}$ at time $(n+1)h_{2}$ conditioned on $\boldsymbol{u}_{n}^{h}=\boldsymbol{r}$ and $\pi_{n}^{h}=c$.
Let $\boldsymbol{U}^{h}$ and $\varPi^{h}$ denote
the collection of ordinary controls, which are determined by measurable
function $\varLambda_{n}^{h}$, such that $\left(\boldsymbol{u}_{n}^{h},\pi_{n}^{h}\right)=\varLambda_{n}^{h}\left(\boldsymbol{\xi}_{k}^{h},k\leq n,\boldsymbol{u}_{k}^{h},\pi_{k}^{h},k<n\right)$.
We say that $\boldsymbol{u}_{n}^{h}$ and $\pi_{n}^{h}$
are admissible for the chain if $\boldsymbol{u}_{n}^{h}$
and $\pi_{n}^{h}$ are $\mathbb{R}_{+}^{d}$
and $[\epsilon,M]$ valued random variables, respectively, and the Markov
property still holds, namely,
\begin{equation*}
\setlength{\abovedisplayskip}{3pt}
\setlength{\belowdisplayskip}{3pt}
\pr\left\{ \boldsymbol{\xi}_{n+1}^{h}=\boldsymbol{z}\left|\boldsymbol{\xi}_{k}^{h},\boldsymbol{u}_{k}^{h},\pi_{k}^{h},k\leq n\right.\right\}=\pr\left\{ \boldsymbol{\xi}_{n+1}^{h}=\boldsymbol{z}\left|\boldsymbol{\xi}_{n}^{h},\boldsymbol{u}_{n}^{h},\pi_{n}^{h}\right.\right\} \coloneqq\pr^{h}\left(\left(\boldsymbol{\xi}_{n}^{h},\boldsymbol{z}\right)\left|\boldsymbol{u}_{n}^{h},\pi_{n}^{h}\right.\right).\end{equation*}
By using the Markov chain above, we can approximate the objective
function defined in (\ref{eq:2-12}) by
\begin{equation*}
\setlength{\abovedisplayskip}{3pt}
\setlength{\belowdisplayskip}{3pt}
J(t,\boldsymbol{y};\boldsymbol{u}^{h},\boldsymbol{\pi}^{h})=\ep_{t}\left[F(\eta_{N_{h_{2}}}^{h})\right]-G\left[\ep_{t}(\eta_{N_{h_{2}}}^{h})\right],
\end{equation*}
where $\eta$ is the first component of $\boldsymbol{\xi}$ and $\eta_{N_{h_{2}}}^{h}$ is the terminal value of the discretized wealth process.

The sequence $\{\boldsymbol{\xi}_{n}^{h},n<\infty\}$ is
said to be locally consistent (see, e.g., \cite{KhDp:13}) w.r.t. (\ref{eq:2-10}), if
it satisfies
\[
\setlength{\abovedisplayskip}{3pt}
\setlength{\belowdisplayskip}{3pt}
\begin{aligned}
\ep_{\boldsymbol{y},n}^{h,\boldsymbol{r},c}\triangle\boldsymbol{\xi}_{n}^{h} & =\boldsymbol{f}(t,\boldsymbol{y},\boldsymbol{r},c)h_{2}+o(h_{2}),\\
\var_{\boldsymbol{y},n}^{h,\boldsymbol{r},c}\triangle\boldsymbol{\xi}_{n}^{h} & =\boldsymbol{\varSigma}(t,\boldsymbol{y},\boldsymbol{r},c)\boldsymbol{\varSigma}(t,\boldsymbol{y},\boldsymbol{r},c)^{\top}h_{2}+o(h_{2}),\\
\underset{n}{\sup}\left|\triangle\boldsymbol{\xi}_{n}^{h}\right| & \rightarrow0,\quad\text{as}\;h\rightarrow0,
\end{aligned}
\]
where $\triangle\boldsymbol{\xi}_{n}^{h}\coloneqq\boldsymbol{\xi}_{n+1}^{h}-\boldsymbol{\xi}_{n}^{h}$,
$\ep_{\boldsymbol{y},n}^{h,\boldsymbol{r},c}$ and $\var_{\boldsymbol{y},n}^{h,\boldsymbol{r},c}$
denote the conditional expectation and variance given by $\{\boldsymbol{\xi}_{k}^{h},\boldsymbol{u}_{k}^{h},\pi_{k}^{h},k\leq n,\boldsymbol{\xi}_{n}^{h}=\boldsymbol{y},\boldsymbol{u}_{n}^{h}=\boldsymbol{r},\pi_{n}^{h}=c\}$,
respectively.

Now, we need to find transition probabilities such that the approximating
Markov chain constructed above is locally consistent. To proceed,
we first suppose that control space has a unique admissible pair control
$\left(\hat{\boldsymbol{u}},\hat{\pi}\right)$ so that
we can drop $\inf$ in (\ref{eq:2-15}). We discretize (\ref{eq:2-15})
by the following finite difference method using step sizes $h=(h_{1},h_{2})$.

Let $\boldsymbol{V}\coloneqq(V,g)^{\top}$ and denote $\boldsymbol{V}^{h}(t,\boldsymbol{y})$ the solution to the finite difference equation, we have
$\boldsymbol{V}(t,\boldsymbol{y})\rightarrow \boldsymbol{V}^{h}(t,\boldsymbol{y})$.

For the derivative w.r.t. time variable, we use
\begin{equation*}
\setlength{\abovedisplayskip}{3pt}
\setlength{\belowdisplayskip}{3pt}
\frac{\partial V^{h}(t,\boldsymbol{y})}{\partial t}\rightarrow\frac{V^{h}(t+h_{2},\boldsymbol{y})-V^{h}(t,\boldsymbol{y})}{h_{2}}.
\end{equation*}
For the first derivative w.r.t. $x$, we use one-side difference
\[
\setlength{\abovedisplayskip}{3pt}
\setlength{\belowdisplayskip}{3pt}
\begin{aligned}\frac{\partial V^{h}(t,\boldsymbol{y})}{\partial x}\rightarrow & \begin{cases}
\frac{V^{h}(t+h_{2},x+h_{1},\boldsymbol{\varphi})-V^{h}(t+h_{2},x,\boldsymbol{\varphi})}{h_{1}}, & \text{if}\;\bar{b}\left(t,\boldsymbol{y},\hat{\boldsymbol{u}},\hat{\pi}\right)\geq0,\\
\frac{V^{h}(t+h_{2},x,\boldsymbol{\varphi})-V^{h}(t+h_{2},x-h_{1},\boldsymbol{\varphi})}{h_{1}}, & \text{if}\;\bar{b}\left(t,\boldsymbol{y},\hat{\boldsymbol{u}},\hat{\pi}\right)<0.
\end{cases}\end{aligned}
\]
For the second derivative w.r.t. $x$, we have standard difference
\[
\setlength{\abovedisplayskip}{3pt}
\setlength{\belowdisplayskip}{3pt}
\begin{aligned}
\frac{\partial^{2}V^{h}(t,\boldsymbol{y})}{\partial x^{2}}
\rightarrow & \frac{\sum_{l=1}^{2}V^{h}(t+h_{2},x+(-1)^{l}h_{1},\boldsymbol{\varphi})-2V^{h}(t+h_{2},x,\boldsymbol{\varphi})}{h_{1}^{2}}.
\end{aligned}
\]
Similarly, we can approximate the first and second derivatives w.r.t. $\varphi^{i}$ and all the partial derivatives of $g^{h}(t,\boldsymbol{y})$ (see also \cite{TkYg:16}).
To simplify notations, we write $\boldsymbol{V}_{n}^{h}(\boldsymbol{y})$ for $\boldsymbol{V}^{h}(nh_{2},\boldsymbol{y})$. Then, detailed calculation leads to the following iterative formula
\[
\setlength{\abovedisplayskip}{3pt}
\setlength{\belowdisplayskip}{3pt}
\begin{aligned} & \boldsymbol{V}_{n}^{h}(\boldsymbol{y})=\left(\tilde{g}_{n+1}^{h}(\boldsymbol{y}),0\right)^{\top}+\boldsymbol{V}_{n+1}^{h}(\boldsymbol{y})p_{1}+\boldsymbol{V}_{n+1}^{h}(x+h_{1},\boldsymbol{\varphi})p_{2+}\\
 & +\boldsymbol{V}_{n+1}^{h}(x-h_{1},\boldsymbol{\varphi})p_{2-}+\sum_{i=1}^{m-1}\left[\boldsymbol{V}_{n+1}^{h}(x,\boldsymbol{\varphi}+h_{1}\boldsymbol{e}_{i})p_{3+}^{i}+\boldsymbol{V}_{n+1}^{h}(x,\boldsymbol{\varphi}-h_{1}\boldsymbol{e}_{i})p_{3-}^{i}\right]\\
 & +\sum_{l=1}^{2}\sum_{k\neq i}^{m-1}\left[\boldsymbol{V}_{n+1}^{h}\left(x,\boldsymbol{\varphi}+(-1)^{l}h_{1}(\boldsymbol{e}_{i}+\boldsymbol{e}_{k})\right)p_{4+}^{i,k}+\boldsymbol{V}_{n+1}^{h}\left(x,\boldsymbol{\varphi}+(-1)^{l}h_{1}(\boldsymbol{e}_{i}-\boldsymbol{e}_{k})\right)p_{4-}^{i,k}\right]
\end{aligned}
\]
where (supressing variables $t,\boldsymbol{y},\hat{\boldsymbol{u}},\hat{\pi}$
of $b$ and $\boldsymbol{\sigma}$)
\begin{align}
p_{1}\coloneqq & \frac{h_{2}}{2h_{1}^{2}}\left[\hat{\pi}\sum_{i,k=1}^{m-1}\varphi^{i}\varphi^{k}\left|\zeta(i)-\bar{\zeta}(\boldsymbol{\varphi})\right|\left|\zeta(k)-\bar{\zeta}(\boldsymbol{\varphi})\right|-3\hat{\pi}\sum_{i=1}^{m-1}\left(\varphi^{i}\left|\zeta(i)-\bar{\zeta}(\boldsymbol{\varphi})\right|\right)^{2}\right]\nonumber \\
 & -\frac{\left(\left|b\right|+\sum_{i=1}^{m-1}\left|\sum_{j=1}^{m}q^{ji}\varphi^{j}\right|\right)h_{1}+\boldsymbol{\sigma}\boldsymbol{\sigma}^{\top}}{h_{1}^{2}}h_{2}+1,\nonumber \\
p_{2\pm}\coloneqq & \frac{\boldsymbol{\sigma}\boldsymbol{\sigma}^{\top}+2b^{\pm}h_{1}}{2h_{1}^{2}}h_{2},\label{eq:3-2}\\
p_{3\pm}^{i}\coloneqq & \frac{h_{2}}{h_{1}^{2}}\left[\hat{\pi}\sum_{i=1}^{m-1}\left(\varphi^{i}\left|\zeta(i)-\bar{\zeta}(\boldsymbol{\varphi})\right|\right)^{2}+(\sum_{j=1}^{m}q^{ji}\varphi^{j})^{\pm}h_{1}-\frac{\hat{\pi}\varphi^{i}\left|\zeta(i)-\bar{\zeta}(\boldsymbol{\varphi})\right|\sum_{k=1}^{m-1}\varphi^{k}\left|\zeta(k)-\bar{\zeta}(\boldsymbol{\varphi})\right|}{2}\right],\nonumber \\
p_{4\pm}^{i,k}\coloneqq & \frac{\hat{\pi}\left[\varphi^{i}\left(\zeta(i)-\bar{\zeta}(\boldsymbol{\varphi})\right)\varphi^{k}\left(\zeta(k)-\bar{\zeta}(\boldsymbol{\varphi})\right)\right]^{\pm}}{4h_{1}^{2}}h_{2}.\nonumber 
\end{align}
Here, $x^{+}$ and $x^{-}$ are the positive and negative parts of
the real number $x$, respectively, and $\tilde{g}_{n+1}^{h}(\boldsymbol{y})$
is given by
\[
\setlength{\abovedisplayskip}{3pt}
\setlength{\belowdisplayskip}{3pt}
\begin{aligned}
\tilde{g}_{n+1}^{h}(\boldsymbol{y}) & =\gamma g_{n+1}^{h}(\boldsymbol{y})\left(1-p_{1}-\frac{\sum_{i=1}^{m-1}\left|\sum_{j=1}^{m}q^{ji}\varphi^{j}\right|+\left|b\right|}{h_{1}}h_{2}\right)\\
 & \quad+\gamma\sum_{l=1}^{2}\left\{ g_{n+1}^{h}(x+(-1)^{l}h_{1},\boldsymbol{\varphi})\frac{b^{+}h_{2}-p_{2+}h_{1}}{h_{1}}+\sum_{i=1}^{m-1}g_{n+1}^{h}(x,\boldsymbol{\varphi}+(-1)^{l}h_{1}\boldsymbol{e}_{i})\frac{(\sum_{j=1}^{m}q^{ji}\varphi^{j})^{+}h_{2}-p_{3+}^{i}h_{1}}{h_{1}}\right.\\
 & \quad\left.-\sum_{k\neq i}^{m-1}\left[\boldsymbol{V}_{n+1}^{h}\left(x,\boldsymbol{\varphi}+(-1)^{l}h_{1}(\boldsymbol{e}_{i}+\boldsymbol{e}_{k})\right)p_{4+}^{i,k}+\boldsymbol{V}_{n+1}^{h}\left(x,\boldsymbol{\varphi}+(-1)^{l}h_{1}(\boldsymbol{e}_{i}-\boldsymbol{e}_{k})\right)p_{4-}^{i,k}\right]\right\} .
\end{aligned}
\]
By choosing proper $h$, we can assume that the coefficient $p_{1}$
is in $[0,1]$, which can be explained as the transition probability
from state $\boldsymbol{y}$ at time $nh_{2}$
to state $\boldsymbol{y}$ at time $(n+1)h_{2}$
given $\left(\hat{\boldsymbol{u}},\hat{\pi}\right)$.
Similarly, we can explain the transition probabilities $p_{2\pm}$,
$p_{3\pm}^{i}$, and $p_{4\pm}^{i,k}$.
Similar to \cite{YzYgZq:15}, we can approximate $V\left(t,\boldsymbol{y}\right)$
by using
\begin{equation*}
\setlength{\abovedisplayskip}{3pt}
\setlength{\belowdisplayskip}{3pt}
V\left(t,\boldsymbol{y}\right)=\underset{\boldsymbol{u}^{h}\in\boldsymbol{U}^{h},\boldsymbol{\pi}^{h}\in\varPi^{h}}{\inf}J\left(t,\boldsymbol{y};\boldsymbol{u}^{h},\boldsymbol{\pi}^{h}\right).
\end{equation*}
It follows from \cite{KhDp:13} that the Markov chain $\{\boldsymbol{\xi}_{n}^{h},n<\infty\}$
with transition probabilities defined in (\ref{eq:3-2})
is indeed locally consistent with (\ref{eq:2-10}).

\section{Approximation of optimal controls}

\subsection{Relaxed controls}

Following \cite{KhDp:13}, let us first recall the definition of
relaxed controls. Let $\mathscr{B}(\boldsymbol{U})$ and $\mathscr{B}(\boldsymbol{U}\times[t,T])$
denote the Borel $\sigma$-algebras of $\boldsymbol{U}$ and $\boldsymbol{U}\times[t,T]$,
respectively. An admissible relaxed control or simply a relaxed control
$m$ is a measure on $\mathscr{B}(\boldsymbol{U}\times[t,T])$
such that $m(\boldsymbol{U}\times[t,s])=s-t$ for all $s\in[t,T]$.
Similarly we can define $\tilde{m}$ a measure on $\mathscr{B}(\varPi\times[t,T])$
such that $\tilde{m}(\varPi\times[t,s])=s-t$. Given relaxed controls
$m$ and $\tilde{m}$, there are derivatives $m_{s}$
and $\tilde{m}_{s}$ which are measures on $\mathscr{B}(\boldsymbol{U})$
and $\mathscr{B}(\varPi)$, respectively, such that $m(d\boldsymbol{r},ds)=m_{s}(d\boldsymbol{r})ds$
and $\tilde{m}(dc,ds)=\tilde{m}_{s}(dc)ds$, respectively.

To prove the convergence of the algorithm, 
we 
use relaxed control representation
\cite{KhDp:13}
Let $\boldsymbol{M}=(M_{1},\dots,M_{d})^{\top}$
and $\tilde{M}$ be a vector-valued measure and
a scalar-valued measure, respectively.
To proceed, we impose the following
conditions.
\begin{description}
\item [{(A1)}] $\boldsymbol{M}$ is square integrable and continuous,
each component is orthogonal, and the pairs $(M_{i},M_{j})_{i \neq j}$
and $\tilde{M}$ are strongly orthogonal.
$\tilde{M}$ is also square integrable, continuous and orthogonal (see, e.g., \cite{YzYgZq:15}).
\end{description}
Under Assumption (A1), there exist measure-valued random processes
$m^{i}$ and $\tilde{m}$ such that the quadratic variation
processes satisfy, for $i=1,2,\dots,d$,
\[
\setlength{\abovedisplayskip}{3pt}
\setlength{\belowdisplayskip}{3pt}
\begin{aligned}
\left\langle M_{i}(A,\cdot),M_{j}(B,\cdot)\right\rangle (s) & =\delta_{ij}m^{i}(A\cap B,s),\\
\left\langle \tilde{M}(\tilde{A},\cdot),\tilde{M}(\tilde{B},\cdot)\right\rangle (s) & =\tilde{m}(\tilde{A}\cap\tilde{B},s),
\end{aligned}
\]
where $A,B\in\mathscr{B}(\boldsymbol{U})$, $\tilde{A},\tilde{B}\in\mathscr{B}(\varPi)$
and $\delta_{ij}$ is Kronecker delta.

\begin{description}
\item [{(A2)}] $m^{i}$ does not depend on $i$,
and $m(\boldsymbol{U},s)=s$.
\end{description}
Under Assumptions (A1) and (A2), there are measures $\boldsymbol{M}$
and $\tilde{M}$ with quadratic variation processes $m\id$
and $\tilde{m}$, respectively, where $\id$ is the unit vector.

With the relaxed controls representation, the operator of the controlled
diffusion is given by
\[
\setlength{\abovedisplayskip}{3pt}
\setlength{\belowdisplayskip}{3pt}
\begin{aligned}\mathscr{L}^{m,\tilde{m}}\phi(t,\boldsymbol{y})\coloneqq & \frac{\partial\phi}{\partial t}+\iint\frac{\partial\phi}{\partial x}b(t,\boldsymbol{y},\boldsymbol{r},c)m_{s}(d\boldsymbol{r})\tilde{m}_{s}(dc)+\frac{1}{2}\int\frac{\partial^{2}\phi}{\partial x^{2}}\boldsymbol{\sigma}(t,\boldsymbol{y},\boldsymbol{r})\boldsymbol{\sigma}^{\top}(t,\boldsymbol{y},\boldsymbol{r})m_{s}(d\boldsymbol{r})\\
 & +\sum_{i=1}^{m-1}\frac{\partial\phi}{\partial\varphi^{i}}\sum_{j=1}^{m}q^{ji}\varphi^{j}+\frac{1}{2}\int\sum_{i,k=1}^{m-1}\varphi^{i}\left(\zeta(i)-\bar{\zeta}(\boldsymbol{\varphi})\right)\varphi^{k}\left(\zeta(k)-\bar{\zeta}(\boldsymbol{\varphi})\right)\frac{\partial^{2}\phi}{\partial\varphi^{i}\partial\varphi^{k}}c\tilde{m}_{s}(dc)\\
= & \iint\mathscr{L}^{\boldsymbol{u},\pi}\phi(s,\boldsymbol{y})m_{s}(d\boldsymbol{r})\tilde{m}_{s}(dc).
\end{aligned}
\]
Besides, there are measures $m$, $\tilde{m}$ satisfying
Assumptions (A1) and (A2) such that for each bounded and smooth function
$\tilde{f}(\cdot,\cdot)$,
\begin{equation*}
\setlength{\abovedisplayskip}{3pt}
\setlength{\belowdisplayskip}{3pt}
\tilde{f}\left(s,\boldsymbol{Y}\right)-\tilde{f}\left(t,\boldsymbol{y}\right)-\iiint\mathscr{L}^{\boldsymbol{u},\pi}\tilde{f}\left(z,\boldsymbol{Y}(z)\right)m_{z}(d\boldsymbol{r})\tilde{m}_{z}(dc)dz
\end{equation*}
is an $\mathbb{\tilde{F}}$ martingale, where $\mathbb{\tilde{F}}=\sigma\left\{ \boldsymbol{Y}(z),m_{z},\tilde{m}_{z},t\leq z\leq s\right\} $.

Then $\left(\boldsymbol{Y},m,\tilde{m}\right)$
solves the martingale problem with operator $\mathscr{L}^{\boldsymbol{u},\pi}$
and we represent our control system as
\begin{equation}
\setlength{\abovedisplayskip}{3pt}
\setlength{\belowdisplayskip}{3pt}
\begin{aligned}
\boldsymbol{Y}(s) & =\boldsymbol{y}+\int_{t}^{s}\int_{\boldsymbol{U}}\int_{\varPi}\boldsymbol{f}\left(\boldsymbol{Y}(z),\boldsymbol{r},c\right)m_{z}(d\boldsymbol{r})\tilde{m}_{z}(dc)dz\label{eq:4-1-1}+\int_{t}^{s}\int_{\boldsymbol{U}}\int_{\varPi}\boldsymbol{\varSigma}\left(\boldsymbol{Y}(z),\boldsymbol{r},c\right)\boldsymbol{M}(d\boldsymbol{r},dz)\tilde{M}(dc,dz).
\end{aligned}
\end{equation}
We say that $\left(\boldsymbol{M},m\right)$ and $\left(\tilde{M},\tilde{m}\right)$
are admissible relaxed controls for (\ref{eq:4-1-1}) if Assumptions
(A1) and (A2) hold and $\left\langle \boldsymbol{M}\right\rangle =m\id$,
$\left\langle \tilde{M}\right\rangle =\tilde{m}$. To
proceed, we pose another two conditions.
\begin{description}
\item [{(A3)}] $\boldsymbol{f}(\cdot,\cdot,\cdot)$ and $\boldsymbol{\varSigma}(\cdot,\cdot,\cdot)$
are continuous; $\boldsymbol{f}(\cdot,\boldsymbol{r},c)$ and $\boldsymbol{\varSigma}(\cdot,\boldsymbol{r},c)$
are Lipschitz continuous uniformly in $\boldsymbol{r},c$ and bounded.
\item [{(A4)}] Let $\boldsymbol{a}\coloneqq\boldsymbol{\varSigma}\boldsymbol{\varSigma}^{\top}$,
there exist a positive constant $C$ and a identity matrix $\boldsymbol{I}$
such that $\boldsymbol{a}-C\boldsymbol{I}$ is positive definite and
$\boldsymbol{a}_{ii}-\sum_{j:j\neq i}\left|\boldsymbol{a}_{ij}\right|\geq0$,
for each $i=1,\dots,m$.
\end{description}

\subsection{Approximation of relaxed controls and measures}
Let $\ep_{n}^{h}$ denote the conditional expectation given
$\{\boldsymbol{\xi}_{k}^{h},\boldsymbol{u}_{k}^{h},\pi_{k}^{h},k\leq n\}$
and define $\boldsymbol{R}_{n}^{h}=\triangle\boldsymbol{\xi}_{n}^{h}-\ep_{n}^{h}\triangle\boldsymbol{\xi}_{n}^{h}$.
By local consistency, we have
\begin{equation*}
\setlength{\abovedisplayskip}{3pt}
\setlength{\belowdisplayskip}{3pt}
\boldsymbol{\xi}_{n+1}^{h}=\boldsymbol{\xi}_{n}^{h}+\boldsymbol{f}(\boldsymbol{\xi}_{n}^{h},\boldsymbol{u}_{n}^{h},\pi_{n}^{h})h_{2}+\boldsymbol{R}_{n}^{h},
\end{equation*}
with
\[
\setlength{\abovedisplayskip}{3pt}
\setlength{\belowdisplayskip}{3pt}
\begin{aligned}
\cov_{n}^{h}\boldsymbol{R}_{n}^{h} & =\boldsymbol{a}(\boldsymbol{\xi}_{n}^{h},\boldsymbol{u}_{n}^{h},\pi_{n}^{h})=\boldsymbol{\varSigma}(\boldsymbol{\xi}_{n}^{h},\boldsymbol{u}_{n}^{h},\pi_{n}^{h})\boldsymbol{\varSigma}^{\top}(\boldsymbol{\xi}_{n}^{h},\boldsymbol{u}_{n}^{h},\pi_{n}^{h})h_{2}+o(h).
\end{aligned}
\]
As shown in \cite{TkYg:16}, by Assumption (A4), we can decompose
\begin{equation*}
\setlength{\abovedisplayskip}{3pt}
\setlength{\belowdisplayskip}{3pt}
\boldsymbol{a}(\boldsymbol{\xi}_{n}^{h},\boldsymbol{u}_{n}^{h},\pi_{n}^{h})=\boldsymbol{T}_{n}^{h}(\boldsymbol{D}_{n}^{h})^{2}(\boldsymbol{T}_{n}^{h})^{\top},
\end{equation*}
where $\boldsymbol{D}_{n}^{h}=\text{diag}(d^{1},\dots,d^{m})$
and $\boldsymbol{T}_{n}^{h}$ is an orthogonal matrix.
Then we can represent the increment of Brownian motion as $\triangle\boldsymbol{W}_{n}^{h}=(\boldsymbol{D}_{n}^{h})^{-1}(\boldsymbol{T}_{n}^{h})^{\top}\boldsymbol{R}_{n}^{h}$, where $\boldsymbol{R}_{n}^{h}=\boldsymbol{\varSigma}(\boldsymbol{\xi}_{n}^{h},\boldsymbol{u}_{n}^{h},\pi_{n}^{h})\triangle\boldsymbol{W}_{n}^{h}+\boldsymbol{\varepsilon}_{n}^{h}$. Thus we have
\begin{equation*}
\setlength{\abovedisplayskip}{3pt}
\setlength{\belowdisplayskip}{3pt}
\boldsymbol{\xi}_{n+1}^{h} =\boldsymbol{\xi}_{n}^{h}+\boldsymbol{f}(\boldsymbol{\xi}_{n}^{h},\boldsymbol{u}_{n}^{h},\pi_{n}^{h})h_{2}+\boldsymbol{\varSigma}(\boldsymbol{\xi}_{n}^{h},\boldsymbol{u}_{n}^{h},\pi_{n}^{h})\triangle\boldsymbol{W}_{n}^{h}+\boldsymbol{\varepsilon}_{n}^{h}.
\end{equation*}
To focus on the control part, let $\{\boldsymbol{L}_{l}^{h},l\leq k_{h}\}$
and $\{C_{\tilde{l}}^{h},\tilde{l}\leq\tilde{k}_{h}\}$
be a finite partition of $\boldsymbol{U}$ and $\varPi$, respectively, such that
$\boldsymbol{L}_{l}^{h}\rightarrow\boldsymbol{0}$ and $C_{\tilde{l}}^{h}\rightarrow0$,
 as $h\rightarrow0$. Let $\boldsymbol{r}_{l}\in \boldsymbol{L}_{l}^{h}$
and $c_{\tilde{l}}\in C_{\tilde{l}}^{h}$, we can define
the random variable
\begin{equation*}
\setlength{\abovedisplayskip}{3pt}
\setlength{\belowdisplayskip}{3pt}
\triangle\boldsymbol{W}_{l,\tilde{l},n}^{h} =\triangle\boldsymbol{W}_{n}^{h}I_{\{\boldsymbol{u}_{n}^{h}=\boldsymbol{r}_{l},\pi_{n}^{h}=c_{\tilde{l}}\}}+\triangle\boldsymbol{\psi}_{l,\tilde{l},n}^{h}I_{\{\boldsymbol{u}_{n}^{h}\neq\boldsymbol{r}_{l},\pi_{n}^{h}\neq c_{\tilde{l}}\}}.
\end{equation*}
Then we have
\begin{align}
 & \boldsymbol{\xi}_{n+1}^{h}=\boldsymbol{\xi}_{n}^{h}+\boldsymbol{f}(\boldsymbol{\xi}_{n}^{h},\boldsymbol{u}_{n}^{h},\pi_{n}^{h})h_{2}\nonumber +\sum_{l,\tilde{l}}\boldsymbol{\varSigma}(\boldsymbol{\xi}_{n}^{h},\boldsymbol{u}_{n}^{h},\pi_{n}^{h})I_{\{\boldsymbol{u}_{n}^{h}=\boldsymbol{r}_{l},\pi_{n}^{h}=c_{\tilde{l}}\}}\triangle\boldsymbol{W}_{l,\tilde{l},n}^{h}+\boldsymbol{\varepsilon}_{n}^{h},\nonumber \\
 & m_{n}^{h}(\boldsymbol{r}_{l})=I_{\{\boldsymbol{u}_{n}^{h}=\boldsymbol{r}_{l}\}},\quad\tilde{m}_{n}^{h}(c_{\tilde{l}})=I_{\{\pi_{n}^{h}=c_{\tilde{l}}\}}.\label{eq:4-2-1}
\end{align}
To approximate the continuous-time following processes $\left(\boldsymbol{Y},m,\boldsymbol{M},\tilde{m},\tilde{M}\right)$, for $s\in\left[nh_{2},(n+1)h_{2}\right)$,
we define the piecewise constant interpolations as
\begin{align}
 & \Psi^{h}(s)=\Psi_{n}^{h},\;\text{for}\;\Psi=\boldsymbol{\xi},\boldsymbol{u},\pi,\boldsymbol{\varepsilon},m,\tilde{m},\quad z^{h_{2}}(s)=n,\nonumber\\
 & \bar{\alpha}^{h}(s)=\sum_{i=1}^{m}g(i)\varphi_{n}^{h,i},\;\boldsymbol{W}_{l,\tilde{l}}^{h}(s)=\sum_{k=0}^{z^{h_{2}}(s)-1}\triangle\boldsymbol{W}_{l,\tilde{l},k}^{h}. \label{eq:4-2-2}
\end{align}
Recall that $m_{n}^{h}$
and $\tilde{m}_{n}^{h}$ are a pair admissible relaxed
controls if $m_{n}^{h}(\boldsymbol{U})=1$, $\tilde{m}_{n}^{h}(\varPi)=1$
and
\[
\setlength{\abovedisplayskip}{3pt}
\setlength{\belowdisplayskip}{3pt}
\begin{aligned}
 & \pr\left\{ \boldsymbol{\xi}_{n+1}^{h}=\boldsymbol{z}\left|\boldsymbol{\xi}_{i}^{h},m_{i}^{h},\tilde{m}_{i}^{h},i\leq n\right.\right\}= \iint\pr^{h}(\boldsymbol{\xi}_{n}^{h},\boldsymbol{z}\left|\boldsymbol{r},c\right.)m_{n}^{h}(d\boldsymbol{r})\tilde{m}_{n}^{h}(dc).
\end{aligned}
\]
For $\boldsymbol{r}_{l}\in\boldsymbol{L}_{l}^{h}$, $c_{\tilde{l}}\in C_{\tilde{l}}^{h}$,
$\{\boldsymbol{M}(\boldsymbol{L}_{l}^{h},\cdot),l\leq k_{h}\}$
and $\{\tilde{M}(C_{\tilde{l}}^{h},\cdot),\tilde{l}\leq\tilde{k}_{h}\}$
are orthogonal continuous martingales with $\left\langle \boldsymbol{M}(\boldsymbol{L}_{l}^{h},\cdot)\right\rangle =m\left(\boldsymbol{L}_{l}^{h},\cdot\right)\id$, and $\left\langle \tilde{M}(C_{\tilde{l}}^{h},\cdot)\right\rangle =\tilde{m}\left(C_{\tilde{l}}^{h},\cdot\right)$.
There are $d+1$ dimensional standard Brownian motions $\boldsymbol{W}_{l,\tilde{l}}^{h}(\cdot),l\leq k_{h},\tilde{l}\leq\tilde{k}_{h}$
such that
\begin{equation*}
\setlength{\abovedisplayskip}{3pt}
\setlength{\belowdisplayskip}{3pt}
\boldsymbol{M}\left(\boldsymbol{L}_{l}^{h},s\right)\tilde{M}(C_{\tilde{l}}^{h},s) =\int_{t}^{s}m_{z}^{1/2}(\boldsymbol{L}_{l}^{h})\tilde{m}_{z}^{1/2}(C_{\tilde{l}}^{h})d\boldsymbol{W}_{l,\tilde{l}}^{h}(z).
\end{equation*}
Let $\boldsymbol{M}^{h}$ and $m^{h}$
be the restrictions of the measures of $\boldsymbol{M}$ and
$m$, respectively, on the set $\{\boldsymbol{L}_{l}^{h},l\leq k_{h}\}$,
and let $\tilde{M}^{h}$ and $\tilde{m}^{h}$
be the restrictions of the measures of $\tilde{M}$ and $\tilde{m}$,
respectively, on the set $\{C_{\tilde{l}}^{h},\tilde{l}\leq\tilde{k}_{h}\}$.
Similar to \cite{KhDp:13}, we have the following lemma.
\begin{lem}

\label{lem1} Under Assumptions (A1)---(A4), it holds that
\begin{equation*}
\setlength{\abovedisplayskip}{3pt}
\setlength{\belowdisplayskip}{3pt}
\left(\boldsymbol{\xi}^{h},m^{h},\boldsymbol{M}^{h},\tilde{m}^{h},\tilde{M}^{h}\right) \Rightarrow\left(\boldsymbol{Y},m,\boldsymbol{M},\tilde{m},\tilde{M}\right),
\end{equation*}
where the notation $\Rightarrow$ denotes
weak convergence,
and $J(t,\boldsymbol{y};m^{h},\tilde{m}^{h})\rightarrow J(t,\boldsymbol{y};m,\tilde{m})$.
\end{lem}

Lemma \ref{lem1} shows we can approximate the variables $\left(\boldsymbol{Y},m,\boldsymbol{M},\tilde{m},\tilde{M}\right)$
by using $\left(\boldsymbol{\xi}^{h},m^{h},\boldsymbol{M}^{h},\tilde{m}^{h},\tilde{M}^{h}\right)$
satisfying
\begin{equation}
\boldsymbol{\xi}^{h}(s)=\boldsymbol{y}+\int_{t}^{s}\sum_{l,\tilde{l}}\boldsymbol{f}(\boldsymbol{\xi}^{h}(z),\boldsymbol{r}_{l},c_{\tilde{l}})m_{z}(\boldsymbol{L}_{l}^{h})\tilde{m}_{z}(C_{l}^{h})dz +\int_{t}^{s}\sum_{l,\tilde{l}}\boldsymbol{\varSigma}(\boldsymbol{\xi}^{h}(z),\boldsymbol{r}_{l},c_{\tilde{l}})m_{z}^{1/2}(\boldsymbol{L}_{l}^{h})\tilde{m}_{z}^{1/2}(C_{\tilde{l}}^{h})d\boldsymbol{W}_{l,\tilde{l}}^{h}(z)+\boldsymbol{\varepsilon}^{h}(s).\label{eq:4-2-3}
\end{equation}
Let $\mathbb{F}^{h}$ represent the $\sigma$-algebra that is generated by
\[
\setlength{\abovedisplayskip}{3pt}
\setlength{\belowdisplayskip}{3pt}
\left\{ \boldsymbol{\xi}^{h}(z),m_{z}^{h},\tilde{m}_{z}^{h},\boldsymbol{M}^{h}(z),\tilde{M}^{h}(z),\boldsymbol{W}_{l,\tilde{l}}^{h}(z),1\leq l\leq k_{h},1\leq\tilde{l}\leq\tilde{k}_{h},t\leq z\leq s\right\} .
\]
Denote $\Gamma^{h}$ and $\tilde{\Gamma}^{h}$
the sets of admissible relaxed controls $m^{h}$ and
$\tilde{m}^{h}$ w.r.t. the set $\left\{ \boldsymbol{\xi}^{h},\boldsymbol{W}_{l,\tilde{l}}^{h},l\leq k_{h},\tilde{l}\leq\tilde{k}_{h}\right\}$, respectively such that $m_{s}^{h}$ and $\tilde{m}_{s}^{h}$
are fixed probability measures in the interval $\left[nh_{2},(n+1)h_{2}\right)$.
We then rewrite
the value function as
\begin{equation}
\setlength{\abovedisplayskip}{3pt}
\setlength{\belowdisplayskip}{3pt}
V^{h}\left(t,\boldsymbol{y}\right) = \inf_{m^{h}\in\Gamma^{h},\tilde{m}^{h}\in\tilde{\Gamma}^{h}}J^{h}\left(t,\boldsymbol{y};m^{h},\tilde{m}^{h}\right),\label{eq:4-2-4}
\end{equation}
and $J^{h}(t,\boldsymbol{y};m^{h},\tilde{m}^{h})=\ep_{t}\left[F\left(\eta^{h}(T)\right)\right]-G\left[\ep_{t}\left(\eta^{h}(T)\right)\right]$.

Next, we use the weak convergence methods (see, e.g., \cite{KhDp:13})
to obtain the convergence of the algorithm. Let $(\boldsymbol{\xi}^{h},m^{h},\boldsymbol{M}^{h},\tilde{m}^{h},\tilde{M}^{h})$
be a solution of (\ref{eq:4-2-3}), where $\boldsymbol{M}^{h}$
and $\tilde{M}^{h}$ are measures with
respect to the filtration $\mathbb{F}^{h}$, with quadratic
variation processes $m^{h}\id$ and $\tilde{m}^{h}$,
respectively. Then we have the following two theorems, and left the details in Appendices A and B.
\begin{thm}
\label{thm1}Under Assumptions (A1)---(A4), let the approximating
chain $\left\{ \boldsymbol{\xi}_{n}^{h},n<\infty\right\} $
be constructed with transition probabilities defined in (\ref{eq:3-2}), let $\{\boldsymbol{u}_{n}^{h},n<\infty\}$
and $\{\pi_{n}^{h},n<\infty\}$ be sequences of admissible
controls, $\boldsymbol{\xi}^{h}$ be the continuous
time interpolations defined by (\ref{eq:4-2-2}), $m^{h}$
be the relaxed control representation of $\boldsymbol{u}^{h}$
(continuous time interpolation of $\boldsymbol{u}_{n}^{h}$),
and $\tilde{m}^{h}$ be the relaxed control representation
of $\pi^{h}$ (continuous time interpolation of $\pi_{n}^{h}$).
Then $H^{h}\coloneqq\{\boldsymbol{\xi}^{h},m^{h},\boldsymbol{M}^{h},\tilde{m}^{h},\tilde{M}^{h}\}$
is tight, which has a weakly convergent
subsequence with the limit $H\coloneqq\{\boldsymbol{Y},m,\boldsymbol{M},\tilde{m},\tilde{M}\}$.
\end{thm}
\begin{thm}
\label{thm2}Under Assumptions (A1)---(A4), let $V\left(t,\boldsymbol{y}\right)$ and $V^{h}\left(t,\boldsymbol{y}\right)$ be value functions defined
in (\ref{eq:2-12}) and (\ref{eq:4-2-4}), respectively. Then $V^{h}\left(t,\boldsymbol{y}\right)\rightarrow V\left(t,\boldsymbol{y}\right),\text{as}\;h\rightarrow0$.
\end{thm}

\section{A numerical example}
In this section, we provide an example with $m=2$ to demonstrate
our results. Here, we set $h_{1}=0.2$, $h_{2}=0.001$, $\gamma=0.5$, $N=2000$ and  $\pi\in[0.001,2]$.
We consider the same financial market as \cite{YzYgZq:15}.

\begin{figure}[H]
\begin{center}
\includegraphics[height=10cm]{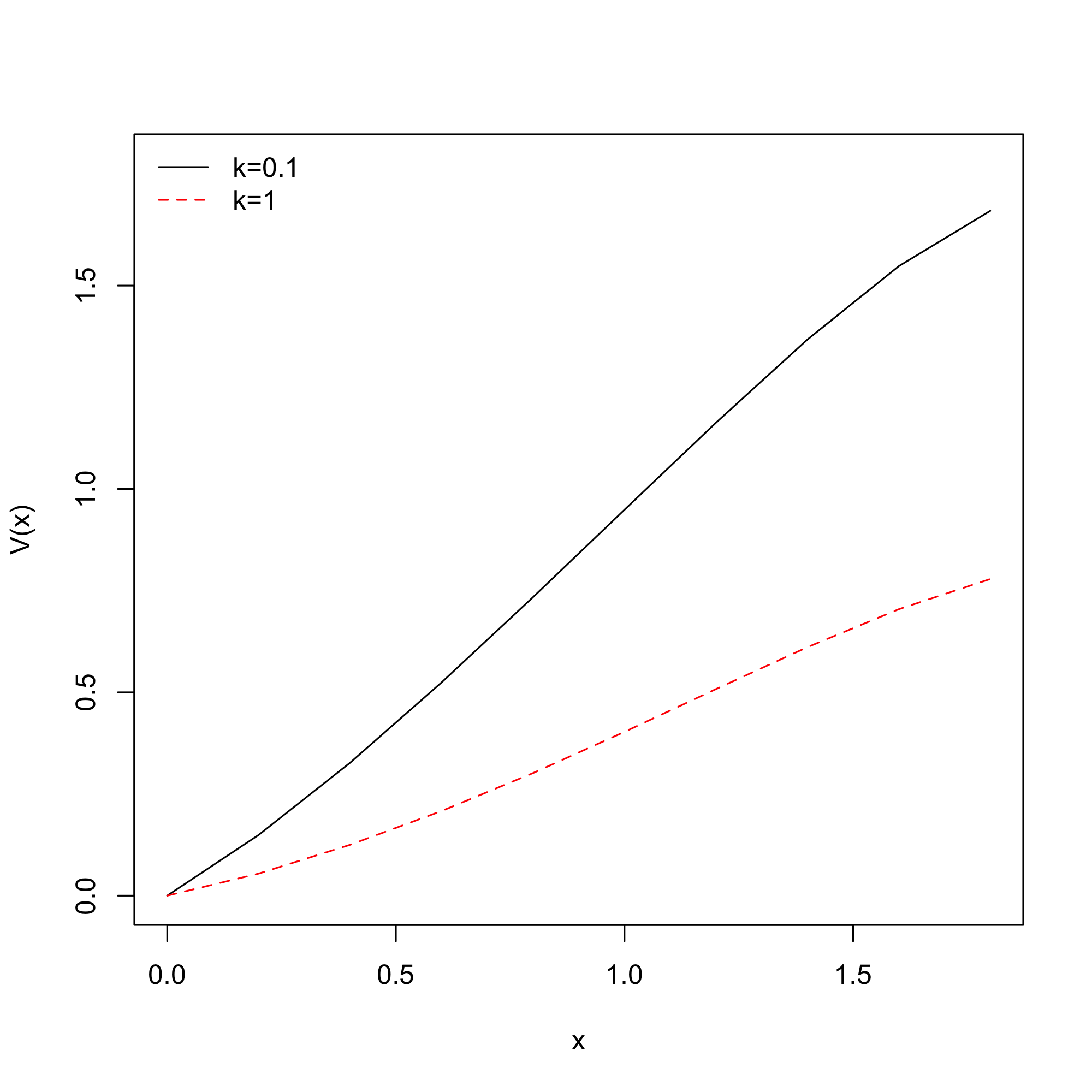}
\caption{Approximate value function $V$.}
\label{fig1}
\end{center}
\end{figure}

\begin{figure}[H]
\begin{center}
\includegraphics[height=10cm]{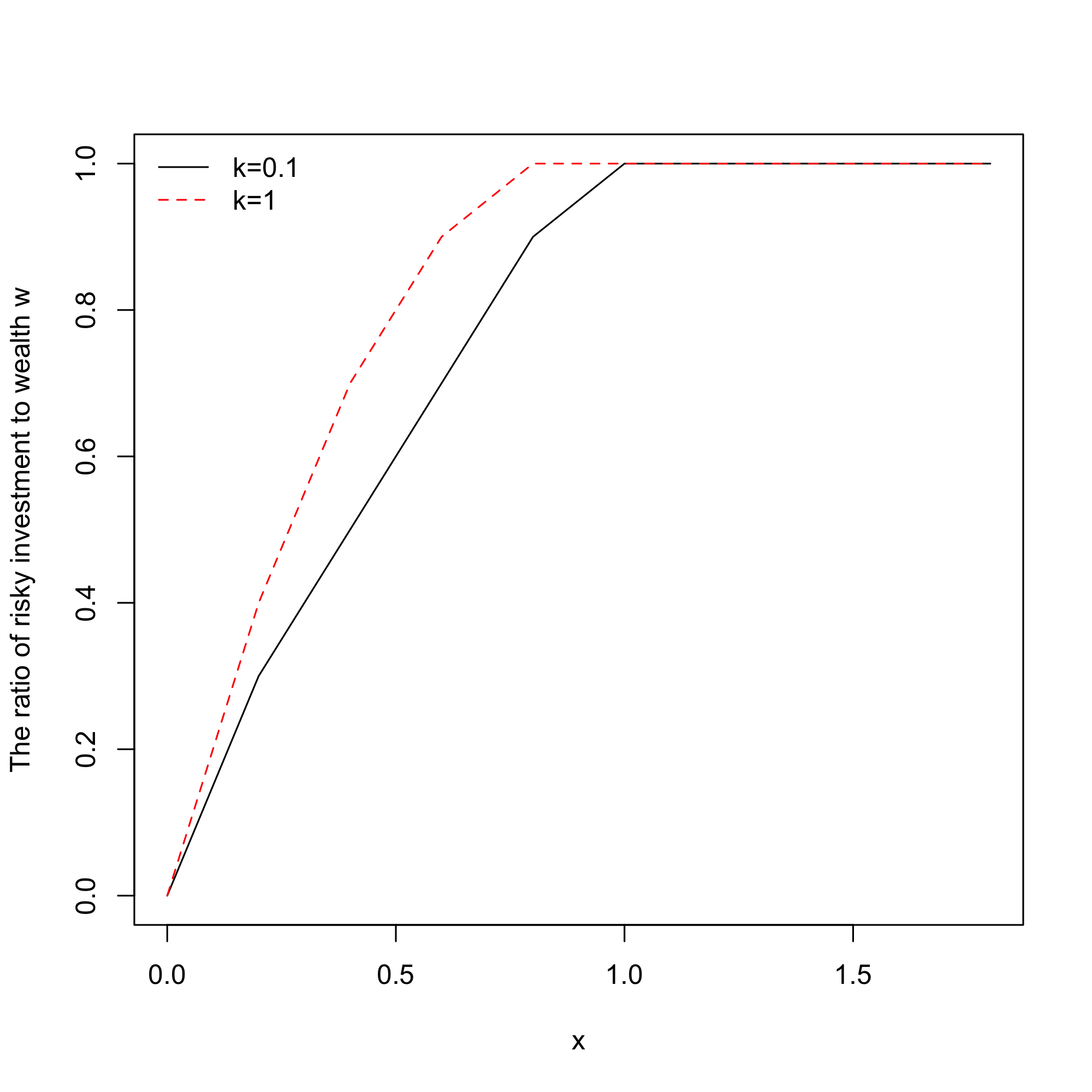}
\caption{The ratio of risky investment to wealth $w$.}
\label{fig2}
\end{center}
\end{figure}

\begin{figure}[H]
\begin{center}
\includegraphics[height=10cm]{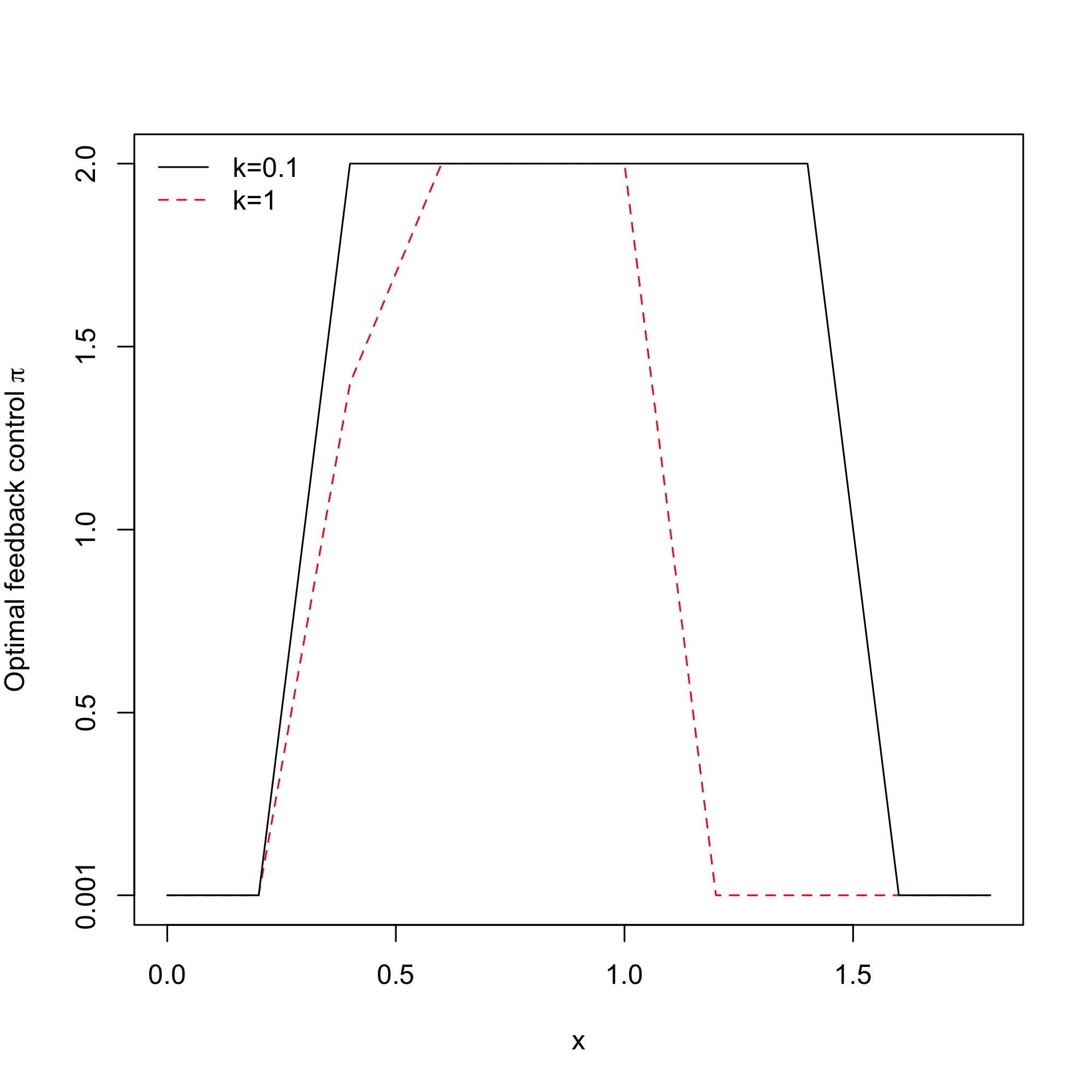}
\caption{Optimal feedback control $\pi$.}
\label{fig3}
\end{center}
\end{figure}

\begin{figure}[H]
\begin{center}
\includegraphics[height=10cm]{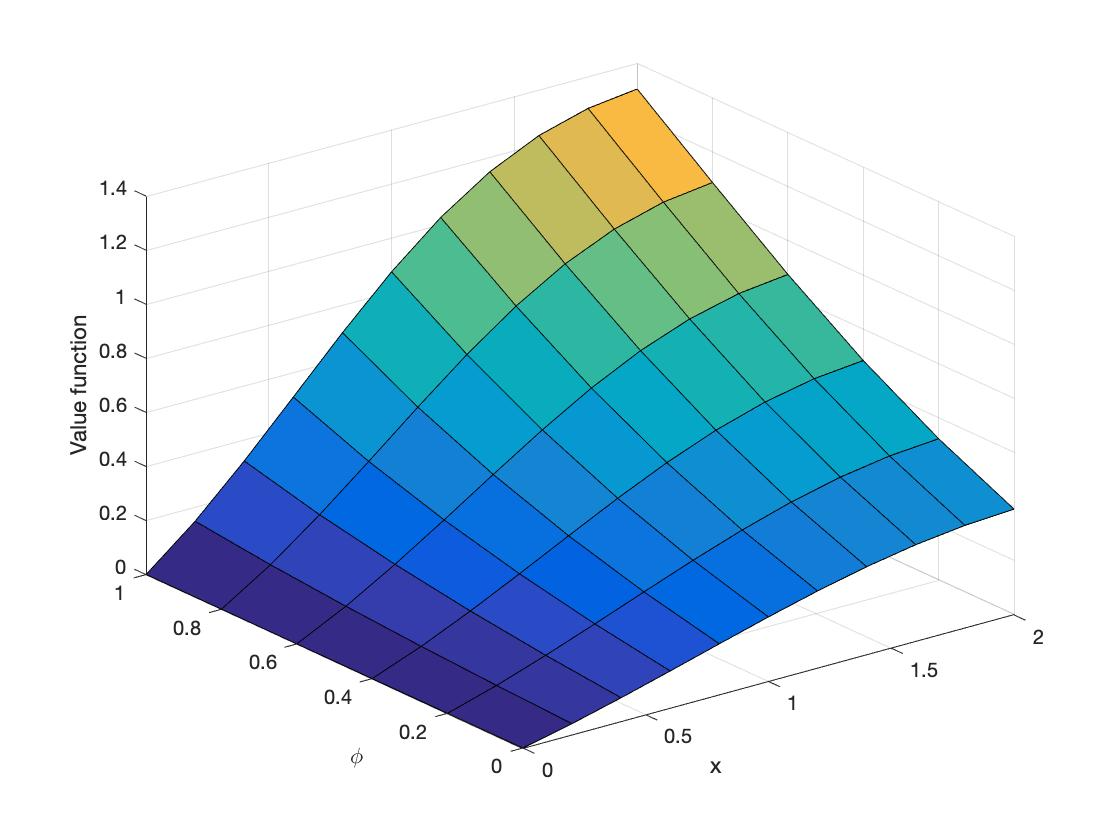}
\caption{Optimal value function versus wealth $x$ and filter state $\varphi$.}
\label{fig4}
\end{center}
\end{figure}

The per-unit-of-wealth information-cost function is specified
in quadratic form, i.e., $\bar{K}(\pi)=k\pi^{2}$, where $k>0$
is the information-cost parameter (see \cite{AdHm:20}).
Fix $t=1$. Let $\phi\coloneqq\varphi^{1}(1)=0.2$, the value function $V(1)$, the ratio of risky
investment to wealth $w(1)$ and the control $\pi(1)$ with different values of $k$ are shown in Figures \ref{fig1}-\ref{fig3}, respectively. The relationship between the value function $V(1)$ and $\left(\phi,x\right)$
is shown in Figure \ref{fig4}.

Figure \ref{fig1} shows the value of value function is larger when
$k$ is smaller.
This is because smaller $k$ implies lower information cost,
and thus the investor can obtain better MV utility.
It follows from Figure \ref{fig2} that an investor intends to invest in risky
assets as wealth increases.
This may be due to the fact
that the investor with less wealth intends to invest in risk-free
assets for risk-averse purposes. Given fixed wealth $x$, larger $k$ implies higher information cost, and thus the investor intends to invest in risky assets, which makes larger $w$.
Figure \ref{fig3} shows that as wealth increases, the investor pays more attention to the signals with small $k$. Because of the specific assumption of the information
cost function in this section, the price of information is relatively
high with larger $k$ when investor's wealth reaches a certain level.
This will lead the investor to reduce investment in the information
market more quickly.

\section{Concluding remarks}
This paper has investigated time-consistent equilibrium strategies
for MV portfolio selection under a hidden Markov model.
We adopted the idea of the dynamic attention behavior in \cite{AdHm:20}
to introduce investor's attention to news on the hidden Markov model.
That is, we considered an investor who can, at each 
time, improve
the accuracy of acquired information at a cost. Under this framework,
we derived an extended HJB equation, for which we used Markov
chain approximation to obtain numerical solutions. 
We constructed
an iterative algorithm, 
proved its convergence, and provided some
numerical results. We also gave some explanations for our numerical
results.


\begin{ack}                               
This work was supported by the Fundamental Research Funds for the Central Universities
 (2020ECNU-HWFW003), the 111 Project (B14019) and the National Natural Science Foundation of China (12071146, 11571113, 11601157, 11601320), Research Grants Council of the Hong Kong Special Administrative Region (project No. 17304921).
 The research of G. Yin was supported  in part by the National Science
Foundation under grant DMS-2204240.
\end{ack}

\bibliographystyle{plain}        


\appendix
\section{Proof of Theorem 4}    
Note that $m^{h}$ and $\tilde{m}^{h}$
are tight due to the compactness of the relaxed controls (see, e.g.,
\cite{KhDp:13}). Here, we consider the tightness of $\boldsymbol{\xi}^{h}$,
by Assumption (A1), for $t\leq s\leq T$,
\[
\begin{aligned}\ep_{t}\left\Vert \boldsymbol{\xi}^{h}(s)-\boldsymbol{Y}(s)\right\Vert ^{2}= & \ep_{t}\left\Vert \int_{t}^{s}\int_{\boldsymbol{U}}\int_{\varPi}\boldsymbol{f}\left(\boldsymbol{\xi}^{h}(z),\boldsymbol{r},c\right)m_{z}^{h}(d\boldsymbol{r})\tilde{m}_{z}^{h}(dc)dz\right.\\
 & \left.+\int_{t}^{s}\int_{\boldsymbol{U}}\int_{\varPi}\boldsymbol{\varSigma}\left(\boldsymbol{\xi}^{h}(z),\boldsymbol{r},c\right)\boldsymbol{M}^{h}(d\boldsymbol{r},dz)\tilde{M}^{h}(dc,dz)+\boldsymbol{\varepsilon}^{h}(s)\right\Vert ^{2}\\
\leq & Ks^{2}+Ks+\varepsilon^{h}(s),
\end{aligned}
\]
where $K$ is a positive constant and for any $t\leq s\leq T$, $\limsup_{h\rightarrow0}\ep\left| \varepsilon^{h}(s)\right| \rightarrow0$. Similarly, we can guarantee $\ep_{t}\left\Vert \boldsymbol{\xi}^{h}(s+\delta)-\boldsymbol{\xi}^{h}(s)\right\Vert ^{2}=O(\delta)+\varepsilon^{h}(\delta)$.
Therefore, the tightness of $\boldsymbol{\xi}^{h}$
follows. By the compactness of sets $\boldsymbol{U}$ and $\varPi$,
we can see that $\boldsymbol{M}^{h}$ and $\tilde{M}^{h}$
are also tight. Hence $H^{h}$ is tight, then by
Prohorov's theorem, $H^{h}$ has a weakly convergent
subsequence with the limit $H$. Next, we show that the limit
is the solution of SDEs driven by $(m,\boldsymbol{M},\tilde{m},\tilde{M})$.

For $\delta>0$ and any process $\nu$ define the process $\nu^{\delta}$
by $\nu^{\delta}(s)=\nu(n\delta)$ for $s\in[n\delta,n\delta+\delta)$.
Then by the tightness of $\boldsymbol{\xi}^{h}$,
(\ref{eq:4-2-3}) can be rewritten as
\[
\begin{aligned}\boldsymbol{\xi}^{h_{1},h_{2}}(s)= & \boldsymbol{y}+\int_{t}^{s}\int_{\boldsymbol{U}}\int_{\varPi}\boldsymbol{f}\left(\boldsymbol{\xi}^{h}(z),\boldsymbol{r},c\right)m_{z}^{h}(d\boldsymbol{r})\tilde{m}_{z}^{h}(dc)dz\\
 & +\int_{t}^{s}\int_{\boldsymbol{U}}\int_{\varPi}\boldsymbol{\varSigma}\left(\boldsymbol{\xi}^{h,\delta}(z),\boldsymbol{r},c\right)\boldsymbol{M}^{h}(d\boldsymbol{r},dz)\tilde{M}^{h}(dc,dz)+\boldsymbol{\varepsilon}^{h,\delta}(s),
\end{aligned}
\]
where $\lim_{\delta\rightarrow0}\limsup_{h\rightarrow0}\ep\left\Vert \boldsymbol{\varepsilon}^{h,\delta}(s)\right\Vert \rightarrow0$.
We further assume that the probability space is chosen as required
by the {Skorohod} representation. Therefore, we can assume the sequence
\begin{equation*}
\left( \boldsymbol{\xi}^{h},m^{h},\boldsymbol{M}^{h},\tilde{m}^{h},\tilde{M}^{h}\right ) \rightarrow\left(\boldsymbol{Y},m,\boldsymbol{M},\tilde{m},\tilde{M}\right)\;\text{w.p.1}
\end{equation*}
with a slight abuse of notation. This leads to that as $h\rightarrow0$,
\[
\begin{aligned}\ep\left\Vert \int_{t}^{s}\int_{\boldsymbol{U}}\int_{\varPi}\boldsymbol{f}\left(\boldsymbol{\xi}^{h}(z),\boldsymbol{r},c\right)m_{z}^{h}(d\boldsymbol{r})\tilde{m}_{z}^{h}(dc)dz-\int_{t}^{s}\int_{\boldsymbol{U}}\int_{\varPi}\boldsymbol{f}\left(\boldsymbol{Y}(z),\boldsymbol{r},c\right)m_{z}^{h}(d\boldsymbol{r})\tilde{m}_{z}^{h}(dc)dz\right\Vert \rightarrow & 0\end{aligned}
\]
uniformly in $s$. Also, recall that $(m^{h},\tilde{m}^{h})\rightarrow(m, \tilde{m})$
in the `compact weak' topology if and only if
\[
\begin{aligned}\int_{t}^{s}\int_{\boldsymbol{U}}\int_{\varPi}\phi(\boldsymbol{r},c,z)m_{z}^{h}(d\boldsymbol{r})\tilde{m}_{z}^{h}(dc)dz\rightarrow & \int_{t}^{s}\int_{\boldsymbol{U}}\int_{\varPi}\phi(\boldsymbol{r},c,z)m_{z}(d\boldsymbol{r})\tilde{m}_{z}(dc)dz,\end{aligned}
\]
for any continuous and bounded function $\phi(\cdot)$ with compact
support. Thus, weak convergence and Skorohod representation imply
that as $h\rightarrow0$,
\begin{align}
\int_{t}^{s}\int_{\boldsymbol{U}}\int_{\varPi}\boldsymbol{f}\left(\boldsymbol{Y}(z),\boldsymbol{r},c\right)m_{z}^{h}(d\boldsymbol{r})\tilde{m}_{z}^{h}(dc)dz\rightarrow & \int_{t}^{s}\int_{\boldsymbol{U}}\int_{\varPi}\boldsymbol{f}\left(\boldsymbol{Y}(z),\boldsymbol{r},c\right)m_{z}(d\boldsymbol{r})\tilde{m}_{z}(dc)dz\label{eq:5-1}
\end{align}
uniformly in $s$ on any bounded interval w.p.1.

Recall that $\boldsymbol{M}^{h}$ and $\tilde{M}^{h}$
are measures with quadratic variation processes $m^{h}\id$
and $\tilde{m}^{h}$, respectively, and $\boldsymbol{\xi}^{h,\delta}$
is a piecewise constant function, following from the probability one
convergence, we have
\[
\begin{aligned}\int_{t}^{s}\int_{\boldsymbol{U}}\int_{\varPi}\boldsymbol{\varSigma}\left(\boldsymbol{\xi}^{h,\delta}(z),\boldsymbol{r},c\right)\boldsymbol{M}^{h_{1},h_{2}}(d\boldsymbol{r},dz)\tilde{M}^{h}(dc,dz)\rightarrow & \int_{t}^{s}\int_{\boldsymbol{U}}\int_{\varPi}\boldsymbol{\varSigma}\left(\boldsymbol{Y}^{\delta}(z),\boldsymbol{r},c\right)\boldsymbol{M}^{h}(d\boldsymbol{r},dz)\tilde{M}^{h}(dc,dz).\end{aligned}
\]
Recall that $(\boldsymbol{M}^{h},\tilde{M}^{h})\rightarrow(\boldsymbol{M},\tilde{M})$
in the `compact weak' topology if and only if
\[
\begin{aligned}\int_{t}^{s}\int_{\boldsymbol{U}}\int_{\varPi}\tilde{\phi}(\boldsymbol{r},c,z)\boldsymbol{M}^{h}(d\boldsymbol{r},dz)\tilde{M}^{h}(dc,dz)\rightarrow & \int_{t}^{s}\int_{\boldsymbol{U}}\int_{\varPi}\tilde{\phi}(\boldsymbol{r},c,z)\boldsymbol{M}(d\boldsymbol{r},dz)\tilde{M}(dc,dz),\;\text{as}\;h\rightarrow0,\end{aligned}
\]
for each bounded and continuous function $\tilde{\phi}$.
Similarly, we have
\begin{align}
\int_{t}^{s}\int_{\boldsymbol{U}}\int_{\varPi}\boldsymbol{\varSigma}\left(\boldsymbol{Y}^{\delta}(z),\boldsymbol{r},c\right)\boldsymbol{M}^{h}(d\boldsymbol{r},dz)\tilde{M}^{h}(dc,dz)\rightarrow & \int_{t}^{s}\int_{\boldsymbol{U}}\int_{\varPi}\boldsymbol{\varSigma}\left(\boldsymbol{Y}(z),\boldsymbol{r},c\right)\boldsymbol{M}(d\boldsymbol{r},dz)\tilde{M}(dc,dz)\label{eq:5-2}
\end{align}
uniformly in $s$ on any bounded interval w.p.1.
Combining Equation
(\ref{eq:5-1}) and (\ref{eq:5-2}), we have
\begin{align}
\boldsymbol{Y}(s)= & \boldsymbol{y}+\int_{t}^{s}\int_{\boldsymbol{U}}\int_{\varPi}\boldsymbol{f}\left(\boldsymbol{Y}(z),\boldsymbol{r},c\right)m_{z}(d\boldsymbol{r})\tilde{m}_{z}(dc)dz \nonumber\\
 & +\int_{t}^{s}\int_{\boldsymbol{U}}\int_{\varPi}\boldsymbol{\varSigma}\left(\boldsymbol{Y}^{\delta}(z),\boldsymbol{r},c\right)\boldsymbol{M}(d\boldsymbol{r},dz)\tilde{M}(dc,dz)+\boldsymbol{\varepsilon}^{\delta}(s),\label{eq:5-3}
\end{align}
where $\lim_{\delta\rightarrow0}\ep\left\Vert \boldsymbol{\varepsilon}^{\delta}(s)\right\Vert =0$.
Taking limit of the above equation (\ref{eq:5-3}) as $\delta\rightarrow0$
yields (\ref{eq:4-1-1}).

\section{Proof of Theorem 5}
We first show that $\lim\inf_{h\rightarrow0}V^{h}(t,\boldsymbol{y})\geq V(t,\boldsymbol{y})$.
Let $(\hat{m}^{h},\hat{\tilde{m}}^{h})$ be a
pair optimal relaxed control for $\boldsymbol{Y}^{h}$
for each $h$. That is,
\[
\begin{aligned}
V^{h}(t,\boldsymbol{y}) & =J^{h}(t,\boldsymbol{y};\hat{m}^{h},\hat{\tilde{m}}^{h})=\inf_{m^{h}\in\Gamma^{h},\tilde{m}^{h}\in\tilde{\Gamma}^{h}}J^{h}(t,\boldsymbol{y};m^{h},\tilde{m}^{h}).
\end{aligned}
\]
Choose a subsequence $\{\tilde{h}\}$ of $\{h\}$
such that
\begin{equation*}
\liminf_{h\rightarrow0}V^{h}\left(t,\boldsymbol{y}\right)=\lim_{\tilde{h}}V^{\tilde{h}}\left(t,\boldsymbol{y}\right)=\lim_{\tilde{h}}J^{\tilde{h}}(t,\boldsymbol{y};\hat{m}^{\tilde{h}},\hat{\tilde{m}}^{\tilde{h}}).
\end{equation*}
Let $\{\boldsymbol{\xi}^{\tilde{h}},\hat{m}^{\tilde{h}},\hat{\boldsymbol{M}}^{\tilde{h}_{1},\tilde{h}_{2}},\hat{\tilde{m}}^{\tilde{h}},\hat{\tilde{M}}^{\tilde{h}}\} $
weakly converge to $\{ \boldsymbol{Y},m,\boldsymbol{M},\tilde{m},\tilde{M}\}$. Otherwise, take a subsequence of $\{\tilde{h}\}$
to assume its weak limit. By Theorem \ref{thm1}, Skorohod representation
and dominance convergence theorem, as $\tilde{h}\rightarrow0$,
we have
\[
\begin{aligned}
 & \ep_{t}\left[F\left(\eta^{\tilde{h}}(T)\right)\right]-G\left[\ep_{t}\left(\eta^{\tilde{h}}(T)\right)\right] \rightarrow\ep_{t}\left[F\left(X(T)\right)\right]-G\left[\ep_{t}\left(X(T)\right)\right].
\end{aligned}
\]
Therefore
\begin{equation*}
J^{h}(t,\boldsymbol{y};\hat{m}^{h},\hat{\tilde{m}}^{h})\rightarrow J(t,\boldsymbol{y};m,\tilde{m})\geq V(t,\boldsymbol{y}).
\end{equation*}
It follows that
\begin{equation*}
\liminf_{h\rightarrow0}V^{h}\left(t,\boldsymbol{y}\right)\geq V(t,\boldsymbol{y}).
\end{equation*}
Next, we show that $\lim\sup_{h\rightarrow0}V^{h}(t,\boldsymbol{y})\leq V(t,\boldsymbol{y})$.

With Lemma \ref{lem1}, given any $\rho>0$, there is a $\delta>0$
that we can approximate $(\boldsymbol{Y},m,\boldsymbol{M},\tilde{m},\tilde{M})$
by  $(\boldsymbol{Y}^{\delta},m^{\delta},\boldsymbol{M}^{\delta},\tilde{m}^{\delta},\tilde{M}^{\delta})$
satisfying
\[
\begin{aligned}
 \boldsymbol{Y}^{\delta}(s)=& \boldsymbol{y}+\int_{t}^{s}\int_{\boldsymbol{U}}\int_{\varPi}\boldsymbol{f}\left(\boldsymbol{Y}^{\delta}(z),\boldsymbol{r},c\right)m_{z}^{\delta}(d\boldsymbol{r})\tilde{m}_{z}^{\delta}(dc)dz\\
 & +\int_{t}^{s}\int_{\boldsymbol{U}}\int_{\varPi}\boldsymbol{\varSigma}\left(\boldsymbol{Y}^{\delta}(z),\boldsymbol{r},c\right)\boldsymbol{M}^{\delta}(d\boldsymbol{r},dz)\tilde{M}^{\delta}(dc,dz),
\end{aligned}
\]
where $m^{\delta}$ and $\tilde{m}^{\delta}$ are piecewise
constants and take finitely many values, $\boldsymbol{M}^{\delta}$
and $\tilde{M}^{\delta}$ are represented in terms of a finite
number of $d$-dimensional and $1$-dimensional Brownian
motions, respectively, and the controls are concentrated on the points
$r_{1},r_{2},\dots,r_{N}$ and $c_{1},c_{2},\dots,c_{N}$, respectively,
for all $s$. Let $\hat{\boldsymbol{u}}^{\rho}$, $\hat{\pi}^{\rho}$
be the optimal controls and $m^{\rho}$, $\tilde{m}^{\rho}$
be corresponding relaxed controls representation, and let $\hat{\boldsymbol{Y}}^{\rho}$
be the associated solution process. Since $\hat{m}^{\rho}$ and
$\hat{\tilde{m}}^{\rho}$ are optimal in the chosen class of
controls, we have
\begin{equation}
J(t,\boldsymbol{y};\hat{m}^{\rho},\hat{\tilde{m}}^{\rho})\leq V(t,\boldsymbol{y})+\frac{\rho}{3}.\label{eq:5-4}
\end{equation}
Note that for each given integer $\iota$, there is a measurable function
$\varLambda_{\iota}^{\rho}$ such that
\begin{equation*}
\left(\hat{\boldsymbol{u}}^{\rho}(s),\hat{\pi}^{\rho}(s)\right)=\varLambda_{\iota}^{\rho}\left(\boldsymbol{W}_{l,\tilde{l}}(\tilde{s}),\tilde{s}\leq\iota\delta,l,\tilde{l}\leq N\right),\;\text{on}\;[\iota\delta,\iota\delta+\delta).
\end{equation*}
We next approximate $\varLambda_{\iota}^{\rho}$ by a function
that depends on the sample of $\left(\boldsymbol{W}_{l,\tilde{l}},l,\tilde{l}\leq N\right)$
at a finite number of time points. Let $\theta<\delta$ such that
$\delta/\theta$ is an integer. Because the $\sigma$-algebra determined
by $\left\{ \boldsymbol{W}_{l,\tilde{l}}(v\theta),v\theta\leq\iota\delta,l,\tilde{l}\leq N\right\} $
increases to the $\sigma$-algebra determined by $\left\{ \boldsymbol{W}_{l,\tilde{l}}(\tilde{s}),\tilde{s}\leq\iota\delta,l,\tilde{l}\leq N\right\} $,
the martingale convergence theorem implies that for each $\delta$ and
$\iota$, there are measurable functions $\varLambda_{\iota}^{\rho,\theta}$,
such that as $\theta\rightarrow0$,
\[
\begin{aligned}
 & \varLambda_{\iota}^{\rho,\theta}\left(\boldsymbol{W}_{l,\tilde{l}}(v\theta),v\theta\leq\iota\delta,l,\tilde{l}\leq N\right)=\left(\boldsymbol{u}_{l}^{\rho,\theta},\pi_{\tilde{l}}^{\rho,\theta}\right)\rightarrow\left(\hat{\boldsymbol{u}}^{\rho}(\iota\delta),\hat{\pi}^{\rho}(\iota\delta)\right)\quad\text{w.p.1}.
\end{aligned}
\]
Here, we select $\varLambda_{\iota}^{\rho,\theta}$ such that
there are $N$ disjoint hyper-rectangles that cover the range of its
arguments and that $\varLambda_{\iota}^{\rho,\theta}$ is constant
on each hyper-rectangle. Let $m^{\rho,\theta}$ and $\tilde{m}^{\rho,\theta}$
denote the relaxed controls representation of the ordinary controls
$\boldsymbol{u}^{\rho,\theta}$ and $\pi^{\rho,\theta}$,
respectively, which take values $\boldsymbol{u}_{l}^{\rho,\theta}$
and $\pi_{\tilde{l}}^{\rho,\theta}$ on $[\iota\delta,\iota\delta+\delta)$,
and let $\boldsymbol{Y}^{\rho,\theta}$ denote the associated
solution. Then for small enough $\theta$, we have
\begin{equation}
J\left(t,\boldsymbol{y};m^{\rho,\theta},\tilde{m}^{\rho,\theta}\right)\leq J\left(t,\boldsymbol{y};\hat{m}^{\rho},\hat{\tilde{m}}^{\rho}\right)+\frac{\rho}{3}.\label{eq:5-5}
\end{equation}
Next, we adapt $\varLambda_{\iota}^{\rho,\theta}$ such that
it can be applied to $\boldsymbol{\xi}_{n}^{h}$. Let
controls $\bar{\boldsymbol{u}}_{n}^{h}$, $\bar{\pi}_{n}^{h}$
to be used for the approximation chain $\boldsymbol{\xi}_{n}^{h}$
defined in (\ref{eq:4-2-1}).

For $\iota=1,2,...$ and $n$ such that $nh_{2}\in[\iota\delta,\iota\delta+\delta)$,
we use the controls which are defined by $\left(\bar{\boldsymbol{u}}_{n}^{h},\bar{\pi}_{n}^{h}\right)=\varLambda_{\iota}^{\rho,\theta}\left(\boldsymbol{W}_{l,\tilde{l}}^{h_{1},h_{2}}(v\theta),v\theta\leq\iota\delta,l,\tilde{l}\leq N\right)$.
Recall that $\bar{m}^{h}$ and $\bar{\tilde{m}}^{h}$
denote the relaxed controls representation of the continuous interpolation
of $\bar{\boldsymbol{u}}_{n}^{h}$ and $\bar{\pi}_{n}^{h}$,
respectively, then
\[
\begin{aligned} & \left(\boldsymbol{\xi}^{h},\bar{m}^{h},\boldsymbol{W}_{l,\tilde{l}}^{h},\bar{\tilde{m}}^{h},\varLambda_{\iota}^{\rho,\theta}\left(\boldsymbol{W}_{l,\tilde{l}}^{h}(v\theta),v\theta\leq\iota\delta,l,\tilde{l}\leq N,\iota=0,1,2,\ldots\right)\right)\\
& \rightarrow \left(\boldsymbol{\xi}^{\rho,\theta},\bar{m}^{\rho,\theta},\boldsymbol{W}_{l,\tilde{l}},\bar{\tilde{m}}^{\rho,\theta},\varLambda_{\iota}^{\rho,\theta}\left(\boldsymbol{W}_{l,\tilde{l}}(v\theta),v\theta\leq\iota\delta,l,\tilde{l}\leq N,\iota=0,1,2,\ldots\right)\right).
\end{aligned}
\]
Thus
\begin{equation}
J(t,\boldsymbol{y};\bar{m}^{h},\bar{\tilde{m}}^{h})\leq J(t,\boldsymbol{y};m^{\rho,\theta},\tilde{m}^{\rho,\theta})+\frac{\rho}{3}.\label{eq:5-6}
\end{equation}
Note that
\begin{equation*}
V^{h}(t,\boldsymbol{y})\leq J(t,\boldsymbol{y};\bar{m}^{h},\bar{\tilde{m}}^{h}).
\end{equation*}
Combing the inequalities (\ref{eq:5-4}), (\ref{eq:5-5}) and (\ref{eq:5-6}),
we have $\lim\sup_{h\rightarrow0}V^{h}(t,\boldsymbol{y})\leq V(t,\boldsymbol{y})$
for the chosen subsequence. According to the tightness of $(\boldsymbol{\xi}^{h},\bar{m}^{h},\bar{\tilde{m}}^{h})$
and arbitrary of $\rho$, we get
\begin{equation*}
\limsup_{h\rightarrow0}V^{h}(t,\boldsymbol{y})\leq V(t,\boldsymbol{y}).
\end{equation*}
\end{document}